\documentclass[reqno]{amsart}
\usepackage{graphicx}
\usepackage{epsfig}
\usepackage{subfigure}
\theoremstyle{proposition}
\newtheorem{proposition}{Proposition}
\theoremstyle{definition}

\theoremstyle{remark}

\newcommand{\g}{\mathcal{G}}
\newcommand{\f}{\mathcal{F}}

\newcommand{\beq}{\begin{equation}}
\newcommand{\eeq}{\end{equation}}

\newcommand{\be}{\begin{enumerate}}
\newcommand{\ee}{\end{enumerate}}
\newcommand{\bi}{\begin{itemize}}
\newcommand{\ei}{\end{itemize}}
\newcommand{\bd}{\begin{description}}
\newcommand{\ed}{\end{description}}

\begin{document}

\title{Call option prices based on Bessel processes}


\author{Ju-Yi Yen$^{(1),(2)}$}
\address{(1) Vanderbilt University, Nashville, Tennessee 37240, USA}
\email{ju-yi.yen@vanderbilt.edu}
\address{ (2) Academia Sinica, Taipei, Taiwan}
\curraddr{}
\thanks{}

\author{Marc Yor$^{(3),(4)}$}
\address{(3) Laboratoire de Probabilit\'{e}s et Mod\`{e}les Al\'{e}atoires,
Universit\'{e} Pierre et Marie Curie, Case Courrier 188, 4, Place
Jussieu, 75252 Paris, Cedex 05,  France}
\address{(4) Institut Universitaire de France} \email{}
\thanks{}


\date{\today}

\dedicatory{}

\begin{abstract}
As a complement to some recent work by Pal and Protter \cite{PP},
we show that the call option prices associated with the Bessel
strict local martingales are integrable over time, and we discuss
the probability densities obtained thus.
\end{abstract}

\maketitle

\section{Introduction: some general remarks}
\subsection{}
Let $(M_t, t\ge 0)$ denote a continuous local martingale, taking
values in $\mathbb{R}_+$. To any $K>0$, we associate the process  $(K-M_t)^+, \ t \ge 0$. It is not difficult to show, after
localizing $M$, that this process $(K-M_t)^+$ is a (bounded)
submartingale, and, as a consequence, the function:
$$m_K^{(+)}(t)=E\big[\left(K-M_t\right)^+\big], \ t\ge 0$$ is increasing, and bounded (by $K$).
The study of such functions, considered (essentially) as
distribution functions, has been the subject of the Bachelier
Course \cite{BYa, BYb}, given by the second author. In particular,
if $M_t\underset{t\rightarrow \infty}{\rightarrow}0$, there is the
formula \beq\label{h_trans}
E\left[F_t\left(K-M_t\right)^+\right]=KE\left[F_t 1_{(\g_K\le
t)}\right] \eeq which is valid for every $F_t\ge 0$, $(\f_t)$
measurable, and $\g_K=\sup\left\{t: M_t=K\right\}$. See also
Madan-Roynette-Yor \cite{MRY,MRY1,MRY2}.

\subsection{}
The present paper is devoted to the study of the functions:
$$m_K^{(-)}(t)=E\big[\left(M_t-K\right)^+\big]=E\big[\left(K-M_t\right)^-\big], \ t\ge 0,$$
which play an important role in option pricing, as $(m_K^{(-)}(t))$
is the European call price with strike $K$, and maturity $t$,
associated with the local martingale $(M_t)$. If $(M_t)$ is a
"true" martingale, then $\{(M_t-K)^+\}$ is a submartingale, hence
$(m_K^{(-)}(t), t\ge 0)$ is increasing. On the other hand, if
$(M_t, t\ge 0)$ is a strict local martingale, that is: a local
martingale, which is not a martingale, then the function
$(m_K^{(-)})$ is not in general increasing, or even monotone.

\subsection{}
The most well-known example of a strict local martingale is
$M_t=1/R_t$, where $(R_t, t\ge 0)$
denotes the $BES(3)$ process, starting from 1, or, by scaling,
equivalently from any $r>0$. Then, the study of $(m_K^{(-)}(t))$ in this particular case has
been undertaken in a remarkable
paper by S. Pal and P. Protter \cite{PP}; the results of which have strongly motivated the present
paper.

In the present paper, we take up again the study of this function $(m_K^{(-)}(t))$ in this particular case; we show that:
$$\int_{0}^{\infty}dt \ m_K^{(-)}(t) < \infty.$$ Hence, up to a multiplicative constant $(m_K^{(-)}(t), \ t \ge 0)$ is a probability density on $\mathbb{R}_+$; we identify the Laplace
transform of this probability, and describe it as the law of a certain random variable defined
uniquely in terms of $BES(3)$ process. This is done thanks to the Doob $h$-transform understanding
of $BES(3)$ (from Brownian motion, killed when hitting 0), combined with general identity~(\ref{h_trans}).
We refer the reader to Section~\ref{S2} for precise statements.
In Section~\ref{S3}, we develop the same kind of study but this time with $M_t=1/R_t^{(\delta-2)}, \ t \ge 0$, where $(R_t, \ t \ge 0)$ denotes the $BES(\delta)$ process, starting from 1. In Section~\ref{S4}, we present the graphs of the corresponding functions $(m_K^{(-)}(t), \ t \ge 0)$.

\subsection{} To summarize, the main point of this work is to use
the interpretation of the generalized Black-Scholes quantities in
terms of last passage times (formula~(\ref{h_trans})) in the
framework of Bessel processes in order to derive fine properties
of the call option process, as a function of maturity, written for
the strict local Bessel martingales.

\section{Some results about $r_K^{(3)}(t)\equiv E_1^{(3)}\big[\big(\frac{1}{X_t}-K\big)^+\big], \ t\ge 0$}\label{S2}
\subsection{}
In this section, we change notation slightly: $(X_t, \ t \ge 0)$ denotes the canonical process on $C(\mathbb{R}_+,\mathbb{R}_+)$, $W_x$ is Wiener measure such that $W_x(X_0=x)=1$, and $P_x^{(3)}$
is the law of the $BES(3)$ process starting from $x$. In fact, we shall only consider $x=1$ (except mentioned otherwise).

\subsection{}\label{S2_2}
Here are our 3 main results concerning the functions  $r_K^{(3)}(t)$.
\begin{proposition}\label{Prop1}
The following holds: \bi \item[(i)] $(t\rightarrow\infty)$:
$r_0^{(3)}(t)\sim\sqrt{\frac{2}{\pi t}};
 \ \ r_K^{(3)}(t)\sim\frac{1}{3\sqrt{2\pi}t^{3/2}K^2}$
\item[(ii)] $(t\rightarrow0)$: $r_K^{(3)}(t)\rightarrow (1-K)^+$;
$r_1^{(3)}(t)\sim\sqrt{\frac{t}{2\pi}}$. \ei
\end{proposition}

An important consequence of Proposition~\ref{Prop1}, (i), is that, for
$K>0$, the function $(r_K^{(3)}(t), t\ge 0)$ is integrable over
$\mathbb{R}_+$, hence it is, up to a multiplicative constant, a
density of probability on $\mathbb{R}_+$. We now describe this probability.
\begin{proposition}\label{Prop2}
\mbox{} \bi \item[(i)] The function $(3K^2r_K^{(3)}(t), t \ge 0)$
is a probability density on $\mathbb{R}_+$. \item[(ii)] It is the
density of
\begin{equation}\label{LambdaK3}
\Lambda_K^{(3)}\stackrel{\rm{(law)}}{=}(g_1-\widetilde{T}_k) +\widetilde{T}_kU  
\end{equation}
where, on the RHS of~(\ref{LambdaK3}), $k=1/K$, the variables
$g_1$, $\widetilde{T}_k$, and $U$ are independent,
$g_1=\sup\{t:X_t=1\}$, $\widetilde{T}_k$ is the size-biased
sampling of $T_k=\inf\{t:X_t=k\}$\footnotemark\footnotetext{That
is: $\widetilde{T}_k$ satisfies
$E[f(\widetilde{T}_k)]=3K^2E[f(T_k)T_k]$ for every $f:
\mathbb{R}_+ \rightarrow \mathbb{R}_+$, Borel.}, with $g_1$ and
$T_k$ defined with respect to $P_0^{(3)}$, and finally $U$ is
uniform on $[0,1]$. \ei
\end{proposition}
As a further description of the law of $\Lambda_K^{(3)}$, we present its Laplace transform.
\begin{proposition}\label{Prop3}
The Laplace transform of $\Lambda_K^{(3)} is$: (we use again
$k=1/K$)
$$E\left[\exp(-\lambda\Lambda_K^{(3)})\right]=\frac{3}{\lambda k^2}\left(e^{-\sqrt{2\lambda}}\right)
\left(\frac{\sinh(\sqrt{2\lambda}k)}{\sqrt{2\lambda}k}-1\right)$$
\end{proposition}

\subsection{Proofs of Propositions~\ref{Prop1}, \ref{Prop2}, \ref{Prop3}}
\subsubsection{}\label{S2_3_1} The main ingredients of these proofs are the following:
\bi
\item[(i)] the Doob $h$-process relationship between Brownian motion and $BES(3)$, which may be written as:
\begin{equation*}\label{P_13}
P_1^{(3)}|_{\f_t} = \left(X_{t\wedge T_0}\right)\cdot W_1|_{\f_t} \end{equation*}
\item[(ii)] the particular instance of formula~(\ref{h_trans}) with $M_t=X_{t\wedge T_0}$, under $W_1$:
\begin{equation*}\label{W_1}
W_1\left(F_t\left(k-X_{t\wedge
T_0}\right)^+\right)=kW_1\left(F_t1_{(\gamma_k\le
t)}\right) \end{equation*}
with $\gamma_k=\sup\{t<T_0:X_t=k\}$.
\item[(iii)] the time reversal result: $(X_{T_0-t}, \ t \le T_0)$ under $W_1$ is distributed as \\ $(X_t, \
    t\le g_1)$ under $P_0^{(3)}$.
\ei

\subsubsection{}
Thanks to the preceding points, we may now obtain interesting
description of $r_K^{(3)}(t)$ in terms of first and last passage
times. In fact, we obtain: \beq\label{r_K3_W1} r_K^{(3)}(t)=
W_1\left(\gamma_k<t<T_0\right) \eeq or, equivalently, from the
time reversal result in (iii): \beq\label{r_K3_P0}
r_K^{(3)}(t)=P_0^{(3)}(g_1>t)-P_k^{(3)}(g_1>t). \eeq
\underline{Proof of~(\ref{r_K3_W1})}:\\
Combining (i) and (ii) above, we obtain:
\begin{eqnarray*}
r_K^{(3)}(t)&=& W_1\left(\left(1-KX_t\right)^+1_{(t<T_0)}\right) \\
&=& KW_1\left(\left(k-X_{t\wedge T_0 }\right)^+1_{(t<T_0)}\right), \ \ \mbox{with } k=\frac{1}{K}\\
&=& KkW_1\left(\gamma_k\le t < T_0\right)\\
&\equiv& W_1\left(\gamma_k<t<T_0\right), \ \ \mbox{which is } (\ref{r_K3_W1}).
\end{eqnarray*}

\subsubsection{}
\noindent\underline{Proof of Proposition~\ref{Prop2}}:\\
\bi \item[a)] We deduce from~(\ref{r_K3_W1}) that:
$$\int_0^\infty dt \ r_K^{(3)}(t)=W_1\left(T_0-\gamma_k\right)=E_0^{(3)}(T_k),$$ by time reversal.
Using the fact that $(R_t^2-3t, t \ge 0)$ is a $P_0^{(3)}$-martingale, we get:
$$3E_0^{(3)}(T_k)=k^2.$$ Hence, the constant $c$ we were seeking is: $c=\frac{3}{k^2}$, and
$\frac{3}{k^2}r_K^{(3)}(t)\equiv 3K^2r_K^{(3)}(t)$ is a probability density on $\mathbb{R}_+$.

\item[b)] In order to identify a random variable $\Lambda_K^{(3)}$ with distribution $3K^2r_K^{(3)}(t)$,
we go back to~(\ref{r_K3_W1}) and we get, for any $f: \mathbb{R}_+\rightarrow\mathbb{R}_+$, Borel:
\begin{eqnarray*}
&&\int_0^\infty dt \ f(t) \ 3K^2r_K^{(3)}(t)=3K^2W_1\left(\int_{\gamma_k}^{T_0}dt \ f(t)\right)\\
&&=3K^2W_1\left[f\left(\gamma_k+\left(T_0-\gamma_k\right)U\right)\left(T_0-\gamma_k\right)\right]\\
&&=3K^2E_0^{(3)}\left[f\left(\left(g_1-T_k\right)+T_kU\right)T_k\right],
\end{eqnarray*}
where $U$ is uniform on $[0,1]$, independent of $T_k$ and $g_1$.
Hence using the notation of $\widetilde{T}_k$ for the size-biased
sampling of $T_k$, we get that:
$$E\left(f\left(\Lambda_K^{(3)}\right)\right)=
E_0^{(3)}\left[f\left(\left(g_1-\widetilde{T}_k\right)+
\widetilde{T}_kU\right)\right].$$ \ei \hspace{5in}$\Box$

\subsubsection{}
\noindent\underline{Proof of Proposition~\ref{Prop3}}:\\
From formula~(\ref{r_K3_W1}) again, this Laplace transform is:
$$3K^2W_1\left(\int_{\gamma_k}^{T_0}dt \ e^{-\lambda t}\right)=\frac{3K^2}{\lambda}W_1\left(
e^{-\lambda \gamma_k}-e^{-\lambda T_0}\right).$$ The result now
follows from:
$$W_1\left(e^{-\lambda T_0}\right)=e^{-\sqrt{2\lambda}}=W_1\left(e^{-\lambda \gamma_k}\right)
\frac{\sqrt{2\lambda}k}{\sinh(\sqrt{2\lambda}k)}.$$
\hspace{4.24in}$\Box$

\subsubsection{}

\underline{Proof of Proposition~\ref{Prop1}}:\\
\bi \item[(i)] As $t\rightarrow \infty$, we have:
\begin{eqnarray*}
r_0^{(3)}(t)&=&E_1^{(3)}\left(\frac{1}{R_t}\right)=\frac{1}{\sqrt{t}} \ E_{(1/\sqrt{t})}^{(3)}
\left(\frac{1}{R_1}\right)\\
&\underset{t\rightarrow \infty}{\sim}&\frac{1}{\sqrt{t}} \ E_0^{(3)}\left(\frac{1}{R_1}\right)=\sqrt{\frac{2}{\pi t}}.
\end{eqnarray*}
The equivalent, for $K>0$, of $r_K^{(3)}(t)$, as $t\rightarrow
\infty$, is simply the particular case: $\delta=3$ of the result
in (i) of Proposition~\ref{Prop4}. \item[(ii)] The first statement
follows from the convergence in $L^1$ of $\frac{1}{X_t}$ to 1 as
$t\rightarrow 0$; For the second statement, we use:
\begin{eqnarray*}
r_1^{(3)}(t)&=&W_1\left(\left(1-X_t\right)^+1_{(t<T_0)}\right)\\
&=&W_1\left(\left(1-X_t\right)^+\right)-W_1\left(\left(1-X_t\right)^+1_{(T_0<t)}\right).
\end{eqnarray*}
The first term is: $W_0(X_t^+)=\sqrt{\frac{t}{2\pi}}$, whereas the
second one may be majorized by
\begin{eqnarray*}
W_1\left(T_0<t\right)&=&
W_0\left(T_1<t\right)=W_0\left(|X_1|>\frac{1}{\sqrt{t}}\right)\\
&=& \sqrt{\frac{2}{\pi}}\int_{1/\sqrt{t}}^{\infty}du \
e^{-(u^2/2)} \le e^{-1/2t}
\end{eqnarray*}
 \ei
(Note that, in fact, it is not necessary to use
Proposition~\ref{Prop1} to obtain Propositions~\ref{Prop2}
and~\ref{Prop3}; however, it is the result of
Proposition~\ref{Prop1}, (i), which led us to study
$\Lambda_K^{(3)}$).

\subsubsection{} A more direct proof of the estimate:
$r_K^{(3)}(t)=\mathrm{O}\left(\frac{1}{t^{3/2}}\right)$, as
$t\rightarrow \infty$. \\ Since this estimate plays quite some
role in our paper, it seems of some interest to look for a simple
proof of it. We note that:
$$r_K^{(3)}(t)\equiv E_1^{(3)}\left(\left(\frac{1}{R_t}-K\right)^+\right)\le E_0^{(3)}\left(\left(\frac{1}{R_t}-K\right)^+\right)$$
since, from the additivity property of squares of Bessel processes
(Shiga-Watanabe \cite{SW}), a Bessel process with dimension
$\delta$, starting from $a>0$, dominates stochastically a Bessel
process with dimension $\delta$, starting from 0. Now we have:
\begin{eqnarray*}
&&E_0^{(3)}\left[\left(\frac{1}{R_t}-K\right)^+\right] \\
&=& E_0^{(3)}\left[\left(\frac{1}{\sqrt{t}{R_1}}-K\right)^+\right] \ \ \ \ \mbox{(by scaling)} \\
&=&c\int_{0}^{1/\sqrt{t}K}dr \ r^2e^{-r^2/2}\left(\frac{1}{\sqrt{t}r}-K\right), \ \ \ \mbox{for some universal constant c}\\
&\le&  \frac{c}{\sqrt{t}}\int_{0}^{1/\sqrt{t}K}dr \ re^{-r^2/2}\\
&= & \frac{c}{\sqrt{t}}\left(1-\exp\left(-\frac{1}{2tK^2}\right)\right)
\underset{t\rightarrow \infty}{\sim}\frac{c}{2t^{3/2}K^2}.
\end{eqnarray*}

\section{Extending the previous results to: $r_K^{(\delta)}(t)\equiv E_1^{(\delta)}\big[\big(\frac{1}{R_t^{\delta-2}}-K\big)^+\big], \ \delta >2$}\label{S3}
\subsection{}
In this section, $(X_t, \ t\ge 0)$ still denotes the canonical process on $C(\mathbb{R}_+,\mathbb{R}_+)$, and we consider $P_x^{(\delta)}$ the law of the $BES(\delta)$ process, starting from $x$ (which, again, will be taken mainly equal to 1).

\subsection{} As a parallel to Section~\ref{S2}, we offer 3 results concerning the function $r_K^{(\delta)}$. We note $\nu=\frac{\delta}{2}-1$
\begin{proposition}\label{Prop4}
The following holds: \bi \item[(i)] $(t\rightarrow\infty)$:
$r_0^{(\delta)}(t)\sim\frac{1}{t^\nu2^\nu\Gamma(1+\nu)};
 \ \ r_K^{(\delta)}(t)\sim\frac{C_K^{(\delta)}}{t^{\nu+1}}, \ \mbox{with } C_K^{(\delta)}=\frac{1}{2^{\nu+1}(\nu+1)\Gamma(\nu)K^{(1/\nu)}} $
\item[(ii)] $(t\rightarrow0)$:
$r_K^{(\delta)}(t)\rightarrow(1-K)^+$; $r_1^{(\delta)}(t)\sim
\sqrt{t}(2\nu) \frac{1}{\sqrt{2\pi}}$   \ei
\end{proposition}

\begin{proposition}\label{Prop5}
\mbox{} \bi \item[(i)] The function $\left(\delta
K^{2/(\delta-2)}r_K^{(\delta)}(t), t\ge 0\right)$ is a probability
density on $\mathbb{R}_+$. \item[(ii)] It is the density of
\beq\label{Lambda_K}
\Lambda_K^{(\delta)}\stackrel{\rm{(law)}}{=}\left(g_1-\widetilde{T}_k\right)+\widetilde{T}_kU
\eeq where, on the RHS of~(\ref{Lambda_K}),
$k=1/K^{1/(\delta-2)}$, the variables $g_1$, $\widetilde{T}_k$,
and $U$ are independent, $g_1=\sup\{t:X_t=1\}$, $\widetilde{T}_k$
is the size-biased sampling of $T_k=\inf\{t: X_t =k\}$, with $g_1$
and $T_k$ defined with respect to $P_0^{(\delta)}$, and finally
$U$ is uniform on $[0,1]$.
 \ei
\end{proposition}
Finally, we present the Laplace transform of $\Lambda_K^{(\delta)}$.

\begin{proposition}\label{Prop6}
The Laplace transform of $\Lambda_K^{(\delta)}$ is:
$$E\left[\exp(-\lambda\Lambda_K^{(\delta)})\right]=\frac{2\delta K^{2/(\delta-2)}}{\lambda}
\left(\mathcal{K}_\nu (\sqrt{2\lambda})\right) \left\{\frac{\nu
I_\nu(k\sqrt{2\lambda})}{k^\nu}-\frac{1}{\Gamma(\nu)}\left(\frac{\sqrt{2\lambda}}{2}\right)^\nu
\right\}$$ where $I_\nu$ and $\mathcal{K}_\nu$ denote the usual
modified Bessel functions, with parameter $\nu$ (we use
$\mathcal{K}_\nu$ instead of $K_\nu$ so that no confusion with the
strike $K$ may occur).
\end{proposition}

\subsection{Proofs of Propositions~\ref{Prop4}, \ref{Prop5}, \ref{Prop6}}
We follow the rationale of the proofs of Propositions~\ref{Prop1}, \ref{Prop2}, \ref{Prop3}, after extending adequately the points (i), (ii) and (iii) in (\ref{S2_3_1})
from the $BES(3)$ process to $BES(\delta)$ process, for $\delta>2$.

\subsubsection{}
Here are these extensions:
\bi
\item[(i)$_\delta$] the Doob $h$-process relationship between $BES(\delta)$ and $BES(4-\delta)$, killed upon hitting 0, is
    $$P_{1}^{(\delta)}|_{\f_t}=\left(X_{t\wedge T_0}\right)^{2\nu}\cdot P_{1}^{(4-\delta)}|_{\f_t},$$
    where $\delta=2(1+\nu)$.
\item[(ii)$_\delta$] the particular instance of formula (\ref{h_trans}) with $M_t=(X_{t\wedge T_0})^{2\nu}$
under $P_1^{(4-\delta)}$:
$$E_1^{(4-\delta)}\left[F_t\left(C^{2\nu}-\left(X_{t\wedge T_0}\right)^{2\nu}\right)^+\right]=C^{2\nu}E_1^{(4-\delta)}\left[F_t1_{(\gamma_{_C}\le t)}\right]$$
with $\gamma_{_C}=\sup\left\{t<T_0: X_t=C\right\}$.
\item[(iii)$_\delta$] the time reversal result:\\
$(X_{(T_0-t)}, \ t\le T_0 )$ under $P_1^{(4-\delta)}$ is distributed as
$(X_t, t\le g_1)$ under $P_0^{(\delta)}$.
\ei

\subsubsection{} The preceding results lead us to the following descriptions of $r_K^{(\delta)}(t)$
in terms of first and last passage times:
\begin{eqnarray*}
r_K^{(\delta)}(t)&=&P_1^{(4-\delta)}(\gamma_k\le t \le T_0) \hspace{1.85in} (3)_\delta\\
&=& P_0^{(\delta)}(g_1>t)-P_k^{(\delta)}(g_1>t)\hspace{1.44in} (4)_\delta
\end{eqnarray*}
where $k=1/K^{1/(\delta-2)}$
\subsubsection{}
\underline{Proof of Proposition~\ref{Prop5}}:\\
It suffices to follow the arguments of the proof of
Proposition~\ref{Prop2}, replacing everywhere dimensions 1 and 3
by dimensions $(4-\delta)$ and $\delta$. \hspace{1.6in}$\Box$

\subsubsection{}
\underline{Proof of Proposition~\ref{Prop6}}:\\
Starting from $(3)_\delta$ we obtain:
\begin{eqnarray*}
&&E\left[\exp \left(-\lambda \Lambda _K^{(\delta)}\right)\right]=\left(\frac{\delta}{k^2}\right)E_1^{(4-\delta)}
\left[\int_{g_k}^{T_0}dt \ e^{-\lambda t}\right]\\
&&= \frac{\delta}{k^2\lambda}\left(E_1^{(4-\delta)}\left(e^{-\lambda g_k}\right)-E_1^{(4-\delta)}\left(e^{-\lambda T_0}\right)\right)\\
&&=\frac{\delta}{k^2\lambda}E_1^{(4-\delta)} \left(e^{-\lambda T_0}\right)\left\{\frac{1}{E_0^{(\delta)}\left(e^{-\lambda T_k}\right)}-1\right\}\\
&&=\frac{\delta}{k^2\lambda}E_0^{(\delta)} \left(e^{-\lambda g_1}\right)\left\{\frac{1}{E_0^{(\delta)}\left(e^{-\lambda T_k}\right)}-1\right\}
\end{eqnarray*}
with the help of the time reversal result (iii)$_\delta$. Now, classical computations of Laplace transforms for first hitting times and last passage times of Bessel processes yield (see, e.g., Kent \cite{Kent}, Pitman-Yor \cite{PY}):
\begin{eqnarray*}
E_0^{(\delta)} \left(e^{-\lambda g_1}\right)&=&\frac{2}{\Gamma(\nu)}\left(\frac{\sqrt{2\lambda}}{2}\right)^\nu \mathcal{K}_\nu\left(\sqrt{2\lambda}\right)\\
E_0^{(\delta)} \left(e^{-\lambda
T_k}\right)&=&\frac{(\sqrt{2\lambda}k)^\nu}{2^\nu
\Gamma(\nu+1)I_\nu(k\sqrt{2\lambda})}
\end{eqnarray*}
where $\mathcal{K}_\nu$ and $I_\nu$ denote the usual modified Bessel functions with index $\nu$. Finally, we have obtained the formula given in Proposition \ref{Prop6}.

\subsubsection{}
\underline{Proof of Proposition~\ref{Prop4}}:\\
For $K=0$, the result follows by scaling, as in dimension 3. \\For
$K>0$, \bi \item[(i)] We shall use formula $(4)_\delta$ (to obtain
the asymptotic result for $r_K^{(\delta)}(t), \ t\rightarrow
\infty$) together with the formula for the distribution of a last
passage time of a transient diffusion (see, e.g., Pitman-Yor
\cite{PY}, as well as Salminen \cite{Sal1,Sal2}):
$$P_x\left(g_y \in dt\right)=-\frac{s'(y)a(y)}{2s(y)}p_t(x,y)dt.$$ In our present case, we have:
$s(y)=-1/y^{(\delta-2)}$, so that: $-\frac{s'(y)}{s(y)}=\frac{(\delta-2)}{y}$; $a(y)=1$. Thus we have: $$P_x^{(\delta)}(g_y \in dt)=\left(\frac{\delta-2}{2y}\right)p_t^{(\delta)}(x,y)dt,$$
so that, going back to the expression $(4)_\delta$ for $r_K^{(\delta)}$, we get:
\beq\label{r_delta2}
r_K^{(\delta)}(t)=\frac{\delta-2}{2}\int_{t}^{\infty}ds \ \left(p_s^{(\delta)}(0,1)-p_s^{(\delta)}(k,1)\right).
\eeq
Next, we shall use the following explicit formulae:  (with $\nu=\frac{\delta-2}{2}$)
\begin{eqnarray*}
p_t^{(\delta)}(x,y)&=&\frac{1}{t}\left(\frac{y}{x}\right)^\nu \exp\left(-\frac{(x^2+y^2)}{2t}\right)I_\nu\left(\frac{xy}{t}\right)\\
p_t^{(\delta)}(0,y)&=&\frac{1}{2^\nu t^{\nu+1}\Gamma(\nu+1)}y^{2\nu+1}\exp\left(-\frac{y^2}{2t}\right)
\end{eqnarray*} (see, e.g., Revuz-Yor \cite{RY}, Chapter XI). Thus:
\begin{eqnarray*}
&&p_s^{(\delta)}(0,1)-p_s^{(\delta)}(k,1)\\ &=&
\frac{1}{2^\nu s^{\nu+1} \Gamma(\nu+1)}\exp\left(-\frac{1}{2s}\right)-\frac{1}{s}\frac{1}{k^\nu}\exp
\left(-\frac{k^2+1}{2s}\right)I_\nu\left(\frac{k}{s}\right)\\
&=& \exp\left(-\frac{1}{2s}\right)\left(\frac{1}{s}\right)\left(\frac{1}{2^\nu s^\nu \Gamma(\nu+1)}-
\frac{1}{k^\nu}\exp\left(-\frac{k^2}{2s}\right)I_\nu\left(\frac{k}{s}\right)\right)\\
&=& \exp\left(-\frac{1}{2s}\right)\left(\frac{1}{s}\right)\left(\frac{1}{k^\nu}\right)
\left(\left(\frac{k}{s}\right)^\nu\frac{1}{2^\nu\Gamma(1+\nu)}-\exp\left(-\frac{k^2}{2s}\right)
I_\nu\left(\frac{k}{s}\right)\right).
\end{eqnarray*}
From (\ref{r_delta2}), it now remains to study the asymptotic, as $t\rightarrow \infty$, of:
$$r_K^{(\delta)}(t)=\frac{\nu}{k^\nu}\int_{t}^{\infty}\frac{ds}{s}\exp\left(-\frac{1}{2s}\right)
\left(\left(\frac{k}{s}\right)^\nu\frac{1}{2^\nu\Gamma(1+\nu)}-\exp\left(-\frac{k^2}{2s}\right)I_\nu\left(
\frac{k}{s}\right)\right)$$
Making the change of variables: $s=\frac{1}{x}$, we get, from formula (\ref{r_delta2}):
\begin{eqnarray*}
r_K^{(\delta)}(t)&=&\frac{\nu}{k^\nu}\int_{0}^{1/t}\frac{dx}{x}\exp\left(-\frac{x}{2}\right)\\
&&\left((kx)^\nu\frac{1}{2^\nu\Gamma(1+\nu)}-\exp\left(-\frac{k^2x}{2}\right)I_\nu(kx)\right).
\end{eqnarray*}
Letting: $y=kx$, we get:
\beq\label{r_delta3}
r_K^{(\delta)}=\frac{\nu}{k^\nu}\int_{0}^{k/t}\frac{dy}{y}e^{-\frac{y}{2k}}
\left(\frac{y^\nu}{2^\nu \Gamma(1+\nu)}-e^{-\frac{ky}{2}}I_\nu(y)\right)
\eeq
Now, we have (see, e.g., Lebedev \cite{Lebedev}):
$$I_\nu(y)=\left(\frac{y}{2}\right)^\nu\left(\frac{1}{\Gamma(1+\nu)}+\sum_{j=1}^{\infty}\left(\frac{y}{2}\right)^{2j}
\frac{1}{j!\Gamma(\nu+j+1)}\right).$$
Consequently:
\begin{eqnarray*}
\frac{y^\nu}{2^\nu\Gamma(\nu+1)}-e^{-\frac{ky}{2}}I_\nu(y)
&=&\left(\frac{y}{2}\right)^\nu\frac{1}{\Gamma(\nu+1)}\left(1-e^{-\frac{ky}{2}}\right)\\
&&-e^{-\frac{ky}{2}}\left(\frac{y}{2}\right)^\nu\left(\sum_{j=1}^{\infty}\left(\frac{y}{2}\right)^{2j}
\frac{1}{j!\Gamma(\nu+j+1)}\right)\\
&\underset{y\rightarrow 0}{\sim}&\left(\frac{y}{2}\right)^\nu\left(\frac{1}{\Gamma(1+\nu)}\left(\frac{ky}{2}\right)\right).
\end{eqnarray*}
Finally, going back to (\ref{r_delta3}), we get:
\begin{eqnarray*}
r_K^{(\delta)}(t)&\sim&\left(\frac{\nu}{k^\nu}\right)\frac{1}{\Gamma(1+\nu)}\left(\frac{k}{2}\right)
\int_{0}^{k/t}dy \left(\frac{y}{2}\right)^\nu\\
&\sim&\frac{1}{t^{\nu+1}}\left(\frac{k^2}{2^{\nu+1}(\nu+1)\Gamma(\nu)}\right).
\end{eqnarray*}
Thus, we have obtained the asymptotic result in (i) of Proposition
\ref{Prop4}, with
$C_K^{(\delta)}=\frac{k^2}{2^{\nu+1}(\nu+1)\Gamma(\nu)}$. Note
that in the particular case $\delta=3$, hence: $\nu=\frac{1}{2}$,
we get: $C_K^{(3)}=\frac{1}{3}\frac{1}{\sqrt{2\pi}}\frac{1}{K^2}$
as claimed in Proposition \ref{Prop1}. \item[(ii)] Now, we study
the asymptotic as $t\rightarrow 0$. The first result is easily
obtained, using the convergence in $L^1$ of
$\frac{1}{R_t^{\delta-2}}$ to 1, as $t\rightarrow 0$. \\We now
give the details in the case $K=1$. As previously, we use:
\beq\label{r_1delta}
r_1^{(\delta)}(t)=E_1^{(4-\delta)}\left(\left(1-R_t^{2\nu}\right)_+1_{(t<T_0)}\right).
\eeq Now, $1-R_{t\wedge T_0}^{2\nu}, \ t\ge 0$, is a martingale
under $P_1^{(4-\delta)}$, which may be written as:
\beq\label{N_tnu} N_t^{(\nu)}=(2\nu)\gamma_{\int_{0}^{t\wedge
T_0}ds\left(R_s^{(2\nu-1)}\right)^2 }\eeq with $(\gamma_u, u\ge
0)$ a standard Brownian motion. We then write:
\begin{equation*}
r_1^{(\delta)}(t)=E_1^{(4-\delta)}\left(\left(N_t^{(\nu)}\right)_+\right)-
E_1^{(4-\delta)}\left(\left(N_t^{(\nu)}\right)_+1_{(T_0<t)}\right)
\end{equation*} We shall show:
\beq\label{E_1}
E_1^{(4-\delta)}\left(\left(N_t^{(\nu)}\right)_+\right)\sim
\sqrt{t}(2\nu)\frac{1}{\sqrt{2\pi}} \eeq as well as:
\beq\label{P_1}
P_1^{(4-\delta)}(T_0<t)=P_0^{(\delta)}(g_1<t)=\mathrm{o}(\sqrt{t}).\eeq
Indeed, from (\ref{N_tnu}), we deduce that the $LHS$ of
(\ref{E_1}) is:
\begin{eqnarray*}
&&(2\nu)\sqrt{t}E\left[\left(\widetilde{\gamma}_{\frac{1}{t}\int_{0}^{t\wedge
T_0}ds\left(R_s^{(2\nu-1)}\right)^2 }\right)_+\right]\\
\sim&&(2\nu)\sqrt{t}E\left[\left(\widetilde{\gamma}_1\right)_+\right]\\
=&&(2\nu)\sqrt{t}\frac{1}{\sqrt{2\pi}}
\end{eqnarray*}
where $(\widetilde{\gamma}_u, \ u \ge 0)$ denotes a standard
Brownian motion.

Concerning (\ref{P_1}), we use Getoor's result \cite{G}
$$g_1\stackrel{\rm{(law)}}{=}1/2\gamma_\nu,$$
where $\gamma_\nu$ denotes a gamma($\nu$) variable; thus:
\begin{eqnarray*}
&&P_0^{(\delta)}(g_1<t)=P\left(\gamma_\nu>\frac{1}{2t}\right)\\
=&&\frac{1}{\Gamma(\nu)}\int_{1/2t}^{\infty}dx \ x^{\nu-1}
e^{-x}\\
=&&\frac{1}{\Gamma(\nu)}e^{-(\frac{1}{2t})}\int_{0}^{\infty}dy \
e^{-y}\left(\frac{1}{2t}+y\right)^{\nu-1}.
\end{eqnarray*}
For $\nu \ge 1$, this quantity is equivalent to:
$\frac{1}{\Gamma(\nu)}e^{-(1/2t)}\left(\frac{1}{2t}\right)^{\nu-1}$
whereas, for $\nu<1$, it is majorized by: $e^{-(1/2t)}$. In any
case, we easily deduce (\ref{P_1}).
 \ei
\underline{Remark}: Using (\ref{r_1delta}) and (\ref{P_1}), we see
that: $$r_1^{(\delta)}(t)\sim
E_1^{(4-\delta)}\left(\left(1-R_t^{2\nu}\right)_+\right).$$ It
should be possible to use the explicit form of the semigroup
density $\mathcal{P}_t^{\delta}(1,r)$, for
$\delta^\prime=4-\delta$, at least when $\delta<4$, to derive
(\ref{E_1}) more directly.

\subsubsection{}A more direct proof of the estimate: $r_K^{(\delta)}(t)=\mathrm{O}\left(\frac{1}{t^{\nu+1}}\right)$, as $t\rightarrow \infty$.\\
Similar to the case of dimension 3, we have, for every dimension $\delta>2$:
$$r_K^{(\delta)}(t)\le E_0^{(\delta)}\left(\left(\frac{1}{R_t}-K\right)^+\right)= E_0^{(\delta)}\left(\left(\frac{1}{\sqrt{t}R_1}-K\right)^+\right)$$ so that:
\begin{eqnarray*}
r_K^{(\delta)}(t) &\le& c\int_{0}^{1/\sqrt{t}K}dr \ r^{\delta-1}e^{-r^2/2}\left(\frac{1}{\sqrt{t}r}
\right)\\
&=& \frac{c}{\sqrt{t}}\int_{0}^{1/\sqrt{t}K}dr \ r^{2\nu}e^{-r^2/2}\\
&\le& \frac{c}{\sqrt{t}}\frac{1}{2\nu+1}\left(\frac{1}{\sqrt{t}K}\right)^{2\nu+1}\\
&=& \frac{c}{(2\nu+1)K^{2\nu+1}}\left(\frac{1}{t^{\nu+1}}\right).
\end{eqnarray*}

\section{Drawing the graphs of $(r_K^{(\delta)}(t), t\ge 0)$}\label{S4}
\subsection{} In order to facilitate the drawing of these graphs, we need to use the simplest possible integral representations of these functions. We shall rely essentially upon formula (\ref{r_delta3}) which was the key of our proof of Proposition \ref{Prop4}:
\begin{equation*}\label{r_delta}
r_K^{(\delta)}(t)=\frac{\nu}{k^\nu}\int_{0}^{k/t}\left(\frac{dy}{y}\right)e^{-y/2k}
\left( \frac{y^\nu}{2^\nu\Gamma(1+\nu)}-e^{-\frac{ky}{2}}I_\nu(y)
\right)
\end{equation*}
Using again the decomposition:
$$I_\nu(y)=\left(\frac{y}{2}\right)^\nu\frac{1}{\Gamma(1+\nu)}+\widetilde{I}_\nu(y)$$

$$\mbox{where:} \ \ \ \ \widetilde{I}_\nu(y)=\left(\frac{y}{2}\right)^\nu\left(\sum_{j=1}^{\infty}\left(\frac{y}{2}\right)^{2j}
\frac{1}{j!\Gamma(\nu+j+1)}\right),$$ we obtain:

\begin{equation*}
r_K^{(\delta)}(t)=r_K^{(\delta)1}(t)-r_K^{(\delta)2}(t), \ \
\mbox{with} \end{equation*}

\begin{equation*} \label{r1r2}
\left\{
\begin{array}{ll}
r_K^{(\delta)1}(t) &= \frac{\nu}{k^\nu}\int_{0}^{k/t}\frac{dy}{y}
\left(\frac{y}{2}\right)^\nu\frac{1}{\Gamma(\nu+1)}e^{-(y/2k)}\left(1-e^{-\frac{ky}{2}}\right) \\
r_K^{(\delta)2}(t)&=
\frac{\nu}{k^\nu}\int_{0}^{k/t}\frac{dy}{y}e^{-\frac{y}{2}\left(k+\frac{1}{k}\right)}\widetilde{I}_\nu(y).
\end{array}
\right.
\end{equation*}

\subsection{}
In the case $\delta=3$, a slightly different approach leads to the following:

\begin{equation}\label{r_3}
r_K^{(3)}(t)=\frac{1}{k}\sqrt{\frac{2}{\pi}}\int_{0}^{1/\sqrt{t}}dy
\left[ e^{-y^2/2}-\frac{1}{2ky^2}\left\{
e^{-\frac{(1-k)^2y^2}{2}}-e^{-\frac{(1+k)^2y^2}{2}}\right\}
\right]
\end{equation}
Thus, we see the importance of the function:

\beq
\varphi(A)\stackrel{\rm{def}}{=}\sqrt{\frac{2}{\pi}}\int_{0}^{A}\frac{dy}{y^2}\left(1-e^{-y^2/2}\right),
\ \ \ A\ge0 \eeq and we note that: \beq \varphi(A)=
\widetilde{N}(A)-\frac{1}{A}\sqrt{\frac{2}{\pi}}\left(1-e^{(-A^2/2)}\right)\label{phi_A}
\eeq where:
$$\widetilde{N}(A)=\sqrt{\frac{2}{\pi}}\int_{0}^{A}dy \ e^{-y^2/2}.$$
Indeed, from formula~(\ref{r_3}), we deduce: \beq\label{r_K3}
r_K^{(3)}(t)=\frac{1}{k}\sqrt{\frac{2}{\pi}}\int_{0}^{1/\sqrt{t}}
dy \ e^{-y^2/2}+\frac{1-k}{2k^2} \
\varphi\left(\frac{1-k}{\sqrt{t}}\right)-\frac{1+k}{2k^2} \
\varphi\left(\frac{1+k}{\sqrt{t}}\right) \eeq

From formula~(\ref{r_K3}), we have:
\begin{eqnarray}
r_K^{(3)}(t)&=&\frac{1}{k}\widetilde{N}\left(\frac{1}{\sqrt{t}}\right)+\frac{(1-k)}{2k^2}
\varphi\left(\frac{1-k}{\sqrt{t}}\right)-\frac{(1+k)}{2k^2}\varphi\left(\frac{1+k}{\sqrt{t}}\right) \label{erf} \\
\mbox{(by~(\ref{phi_A})) } &=&
\frac{1}{k}\widetilde{N}\left(\frac{1}{\sqrt{t}}\right)+
\frac{1}{2k^2}\left(\varphi\left(\frac{1-k}{\sqrt{t}}\right)-\varphi\left(\frac{1+k}{\sqrt{t}}\right)
\right) \nonumber \\ &&\ \ \ \ -
\frac{1}{2k}\left(\varphi\left(\frac{1-k}{\sqrt{t}}\right)+\varphi\left(\frac{1+k}{\sqrt{t}}\right)\right)
\nonumber \\
&=&\frac{1}{k}\left\{\widetilde{N}\left(\frac{1}{\sqrt{t}}\right)-\frac{1}{2}
\left[\widetilde{N}\left(\frac{1-k}{\sqrt{t}}\right)+\widetilde{N}\left(\frac{1+k}{\sqrt{t}}\right)
\right.\right. \nonumber \\ && \ \ \ \left. \left.
-\left[\frac{\sqrt{t}}{(1-k)}\sqrt{\frac{2}{\pi}}\left(1-e^{-\frac{(1-k)^2}{2t}}\right)+
\frac{\sqrt{t}}{(1+k)}\sqrt{\frac{2}{\pi}}\left(1-e^{-\frac{(1+k)^2}{2t}}\right)\right]
\right]
\right\} \nonumber\\
&& \ \ \
+\frac{1}{2k^2}\left\{\widetilde{N}\left(\frac{1-k}{\sqrt{t}}\right)-
\widetilde{N}\left(\frac{1+k}{\sqrt{t}}\right)\right. \nonumber \\
&& \ \ \ \left.
-\frac{\sqrt{t}}{(1-k)}\sqrt{\frac{2}{\pi}}\left(1-e^{-\frac{(1-k)^2}{2t}}\right)+
\frac{\sqrt{t}}{(1+k)}\sqrt{\frac{2}{\pi}}\left(1-e^{-\frac{(1+k)^2}{2t}}\right)\right\}\nonumber\\
&\equiv&\frac{1}{k}\left\{\widetilde{N}\left(\frac{1}{\sqrt{t}}\right)-\frac{1}{2}
\left[\widetilde{N}\left(\frac{1-k}{\sqrt{t}}\right)+
\widetilde{N}\left(\frac{1+k}{\sqrt{t}}\right)\right]\right\}
\nonumber \\ && \ \ \ +
\frac{1}{2k^2}\left[\widetilde{N}\left(\frac{1-k}{\sqrt{t}}\right)
-\widetilde{N}\left(\frac{1+k}{\sqrt{t}}\right)\right]
+\frac{\sqrt{t}}{2k(1-k)}\sqrt{\frac{2}{\pi}}\left(1-e^{-\frac{(1-k)^2}{2t}}\right)
\nonumber \\
&& \ \ \  +\frac{\sqrt{t}}{2k(1+k)}\sqrt{\frac{2}{\pi}}\left(1-e^{-\frac{(1+k)^2}{2t}}\right) 
-\frac{\sqrt{t}}{2k^2(1-k)}\sqrt{\frac{2}{\pi}}\left(1-e^{-\frac{(1-k)^2}{2t}}\right)
\nonumber \\
&& \ \ \
+\frac{\sqrt{t}}{2k^2(1+k)}\sqrt{\frac{2}{\pi}}\left(1-e^{-\frac{(1+k)^2}{2t}}\right)\label{r_K3a}
\end{eqnarray}
Formula~(\ref{r_K3a}) involves 5 terms with $\widetilde{N}$, and 4
terms with "$\exp$". Thus, we write:
$$r_K^{(3)}(t)=\widetilde{r}_k(t)+r_k^{(\exp)}(t),$$
with
\begin{eqnarray} \label{r_tilde}
\widetilde{r}_k(t)&=&\frac{1}{k}\left\{\widetilde{N}\left(\frac{1}{\sqrt{t}}\right)-\frac{1}{2}
\left[\widetilde{N}\left(\frac{1-k}{\sqrt{t}}\right)
+\widetilde{N}\left(\frac{1+k}{\sqrt{t}}\right)\right]\right\} \nonumber \\
&& \ \ \ +
\frac{1}{2k^2}\left[\widetilde{N}\left(\frac{1-k}{\sqrt{t}}\right)
-\widetilde{N}\left(\frac{1+k}{\sqrt{t}}\right)\right]
\end{eqnarray}

\begin{eqnarray}
\widetilde{r}_k^{(\exp)}(t)&=&\sqrt{\frac{t}{2\pi}}
\left\{\left(\frac{1}{k}-\frac{1}{k^2}\right)\frac{1}{(1-k)}\left(1-e^{-\frac{(1-k)^2}{2t}}\right)\right.
\nonumber \\ && \ \ \ \left.
+\left(\frac{1}{k}+\frac{1}{k^2}\right)\frac{1}{(1+k)}\left(1-e^{-\frac{(1+k)^2}{2t}}\right)
\right\} \nonumber \\
&=&
\sqrt{\frac{t}{2\pi}}\left\{-\frac{1}{k}\left(1-e^{-\frac{(1-k)^2}{2t}}\right)
+\frac{1}{k}\left(1-e^{-\frac{(1+k)^2}{2t}}\right)\right\} \nonumber \\
&=&
\sqrt{\frac{t}{2\pi}}\frac{1}{k}\left(\exp\left(-\frac{(1-k)^2}{2t}\right)+
\exp\left(-\frac{(1+k)^2}{2t}\right)\right). \label{r_exp}
\end{eqnarray}

We note again the further simplification of~(\ref{r_exp}):
\begin{equation*}
r_K^{(\exp)}(t)=\sqrt{\frac{2t}{\pi}}\left(\frac{1}{k}\right)\left(\exp\left(-\frac{1+k^2}{2t}\right)\right)
\cosh\left(\frac{k}{t}\right). \ \ \ \ \ \ \ (\ref{r_exp}^\prime)
\end{equation*}

\subsection{}
We now present the graphs of $r_K^{\delta}(t)$ for various
dimensions $\delta$ and strikes $K$. In fact, we have drawn two
kinds of graphs: \bi \item[(a)] Figure 1: for each dimension
$\delta$, we consider $k=1, \frac{1}{2},\frac{1}{3}
...,\frac{1}{10}$, and draw all the graphs together. We use
formula (\ref{erf}) for $\delta =3$, and formula (\ref{r_delta3})
for $\delta=5,7,9,11,13$ to draw the graphs. \item[(b)] Figure 2:
for each $k$, we draw the graphs for different dimensions. Here,
we use formula (\ref{r_delta3}) to draw graphs for all dimensions.
\ei

\begin{figure}[htp]
\centering
\begin{tabular}{cc}
\epsfig{file=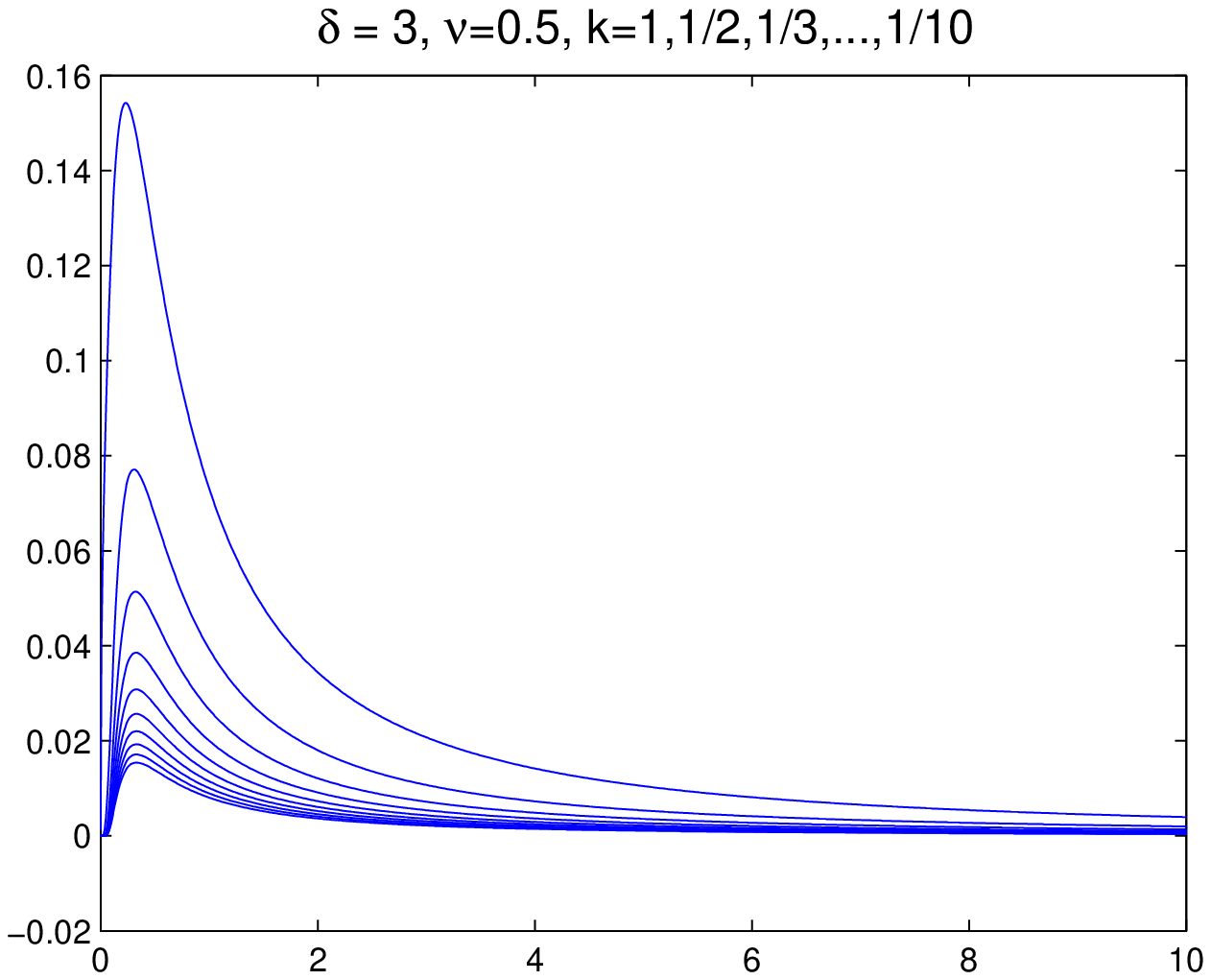,width=0.5\linewidth,clip=} &
\epsfig{file=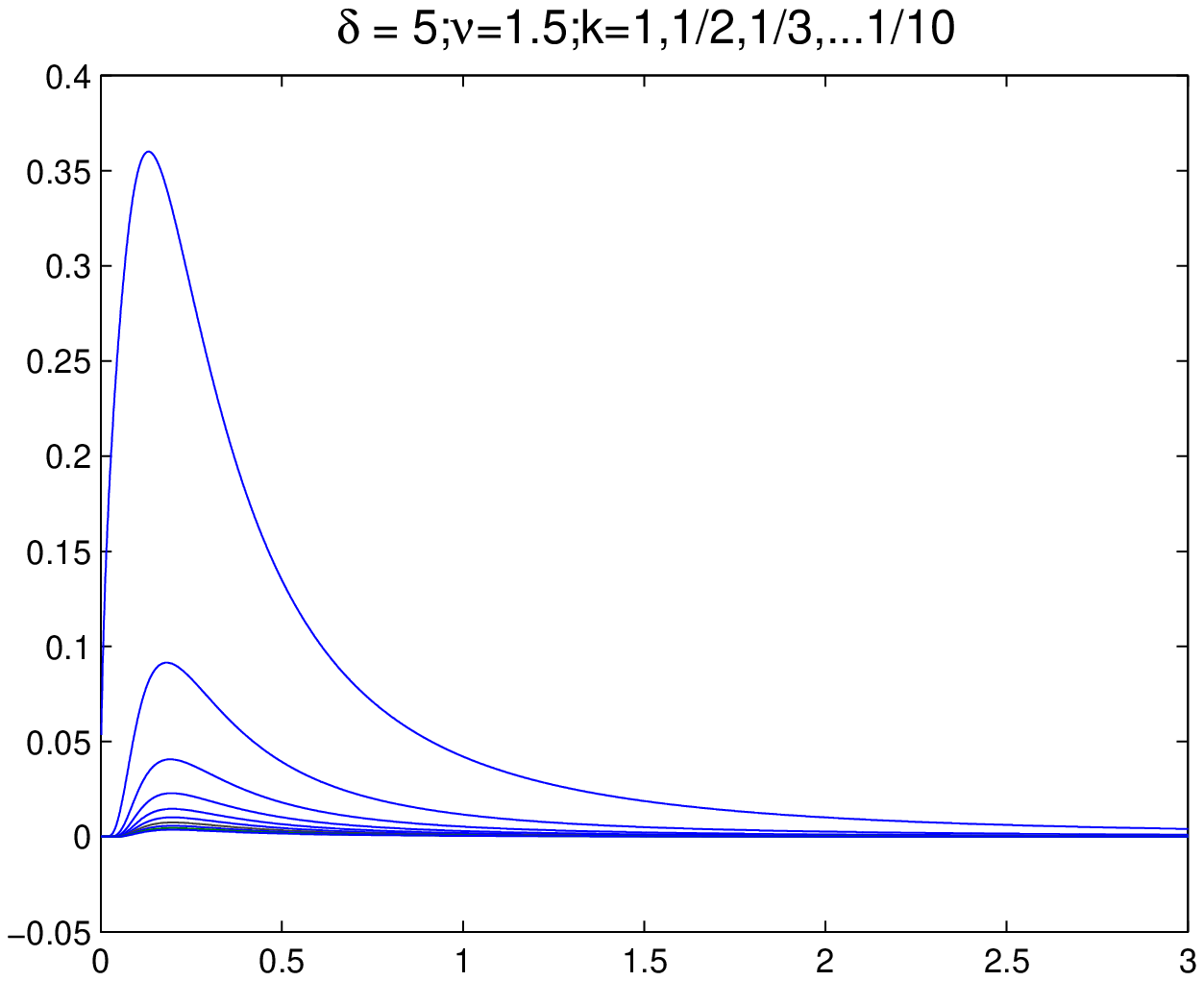,width=0.5\linewidth,clip=} \\
\epsfig{file=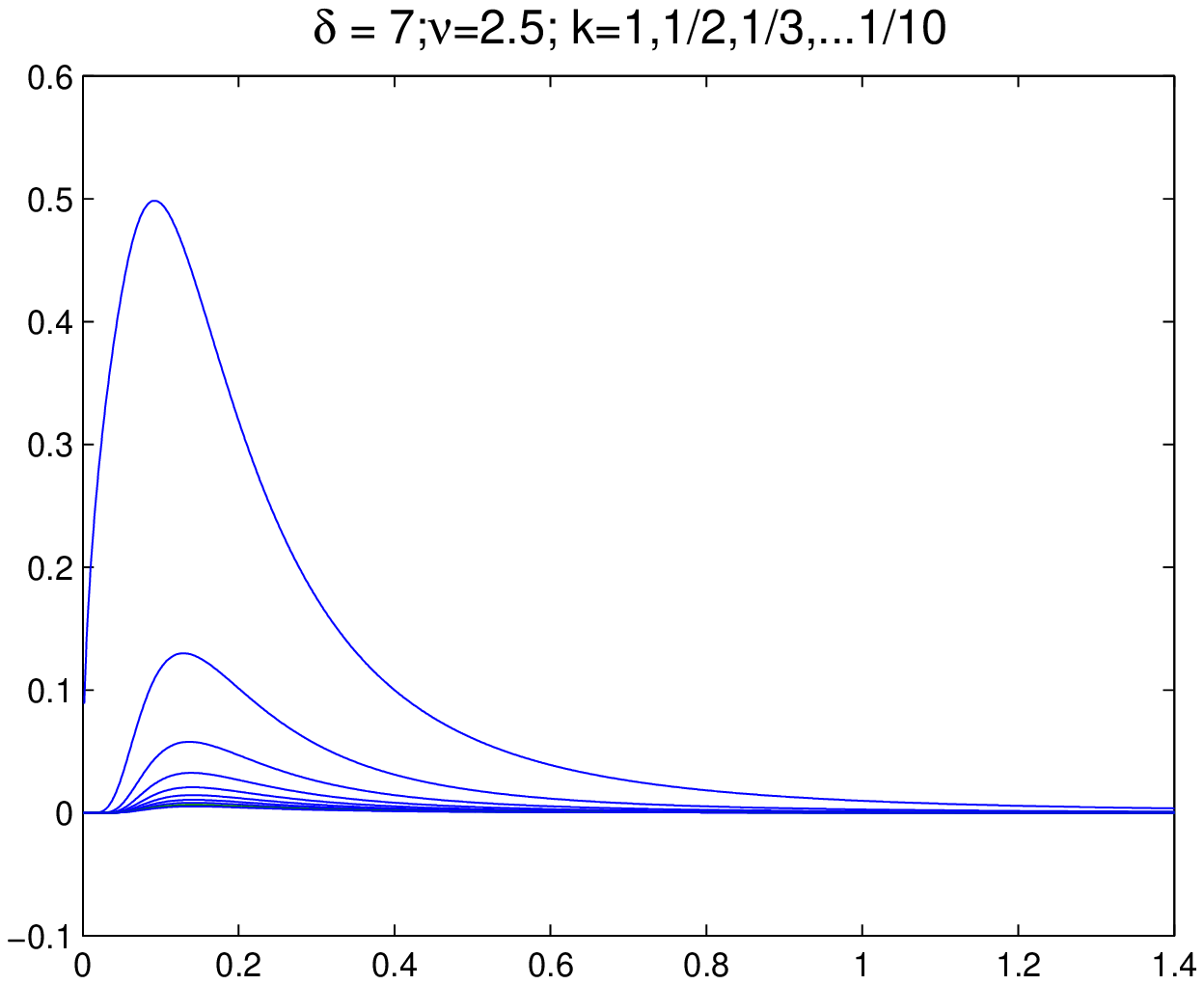,width=0.5\linewidth,clip=} &
\epsfig{file=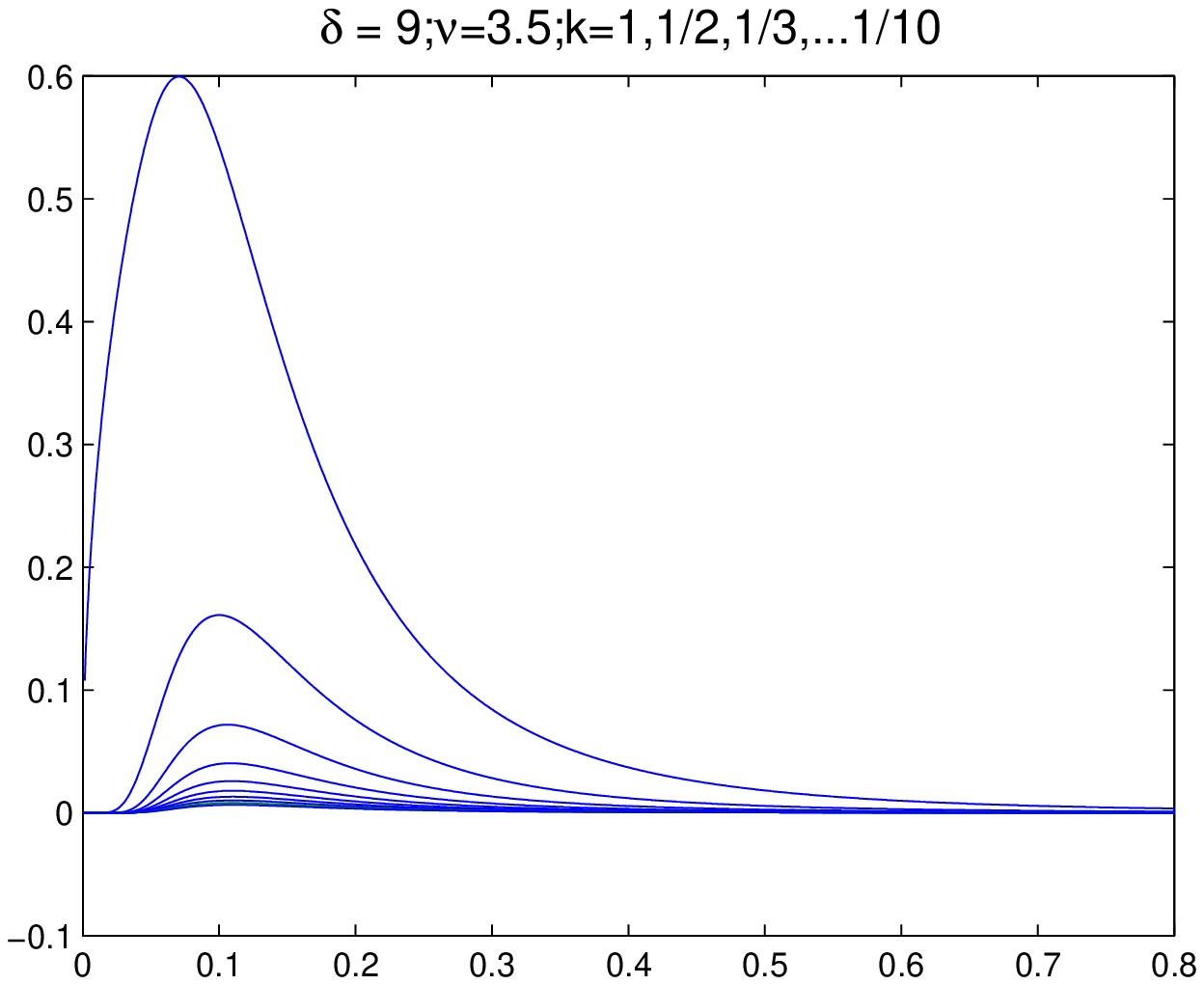,width=0.5\linewidth,clip=}\\
\epsfig{file=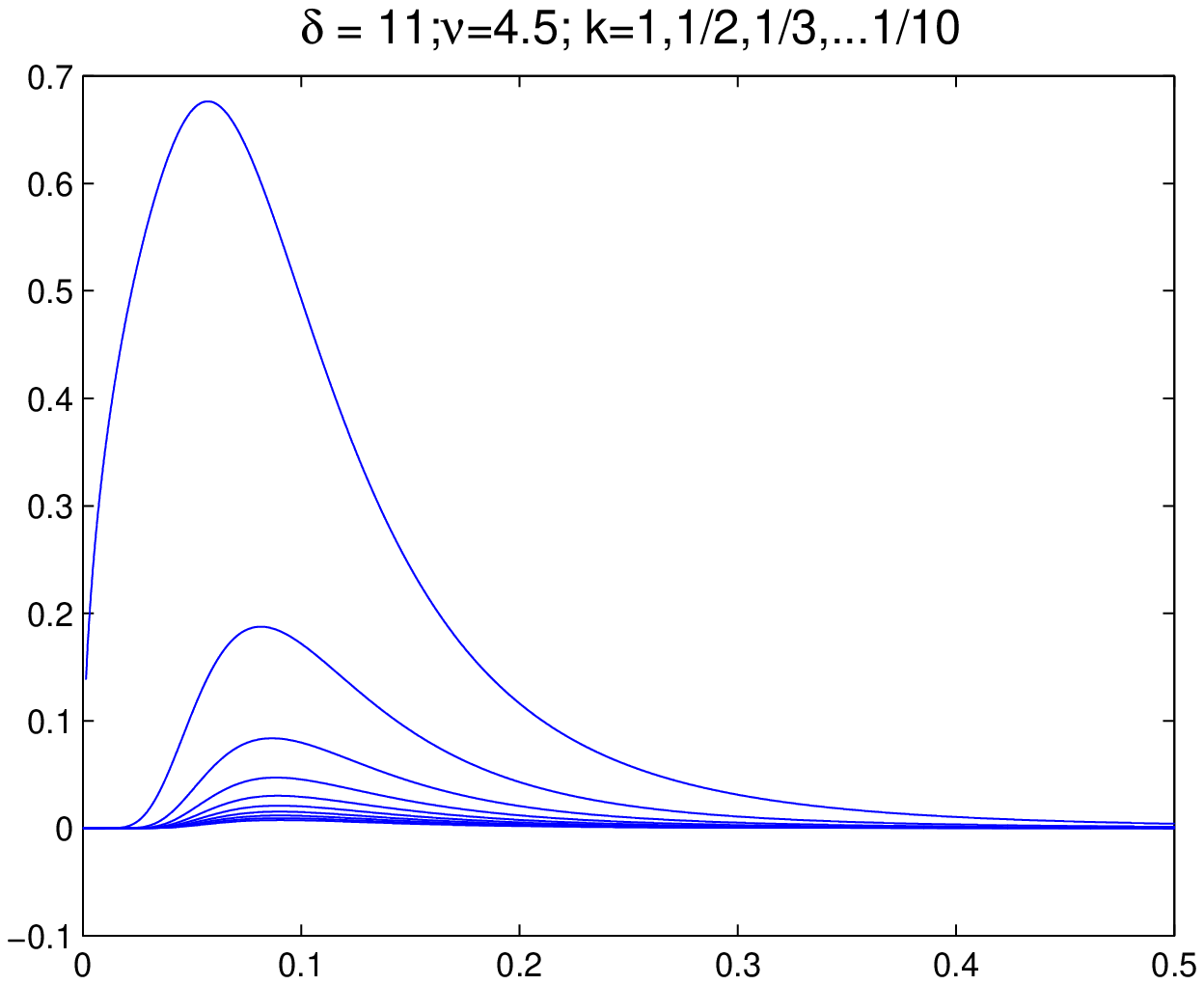,width=0.5\linewidth,clip=} &
\epsfig{file=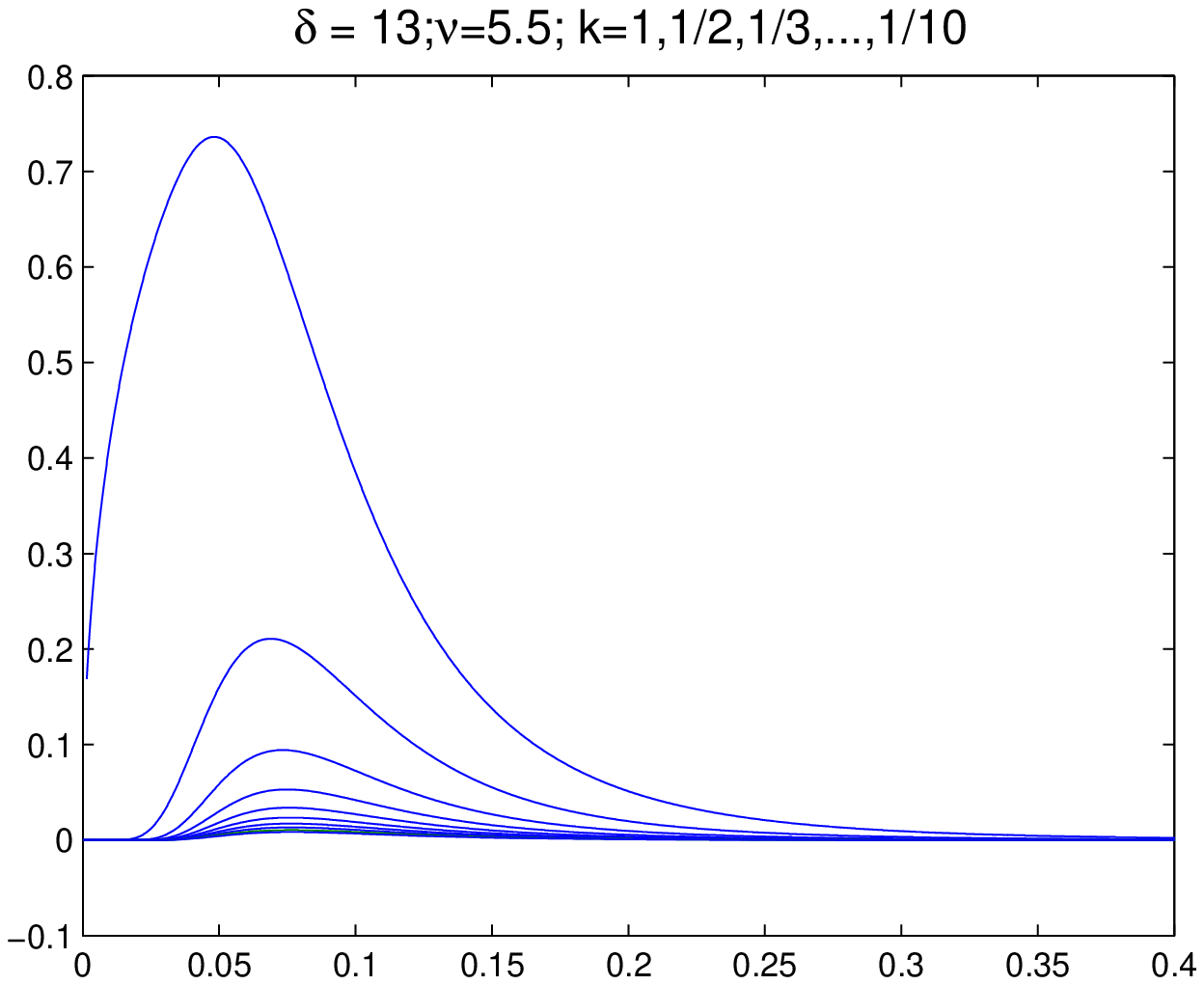,width=0.5\linewidth,clip=}
\end{tabular}
\caption{Varying $k$}
\end{figure}

\begin{figure}[htp]
\centering
\begin{tabular}{cc}
\epsfig{file=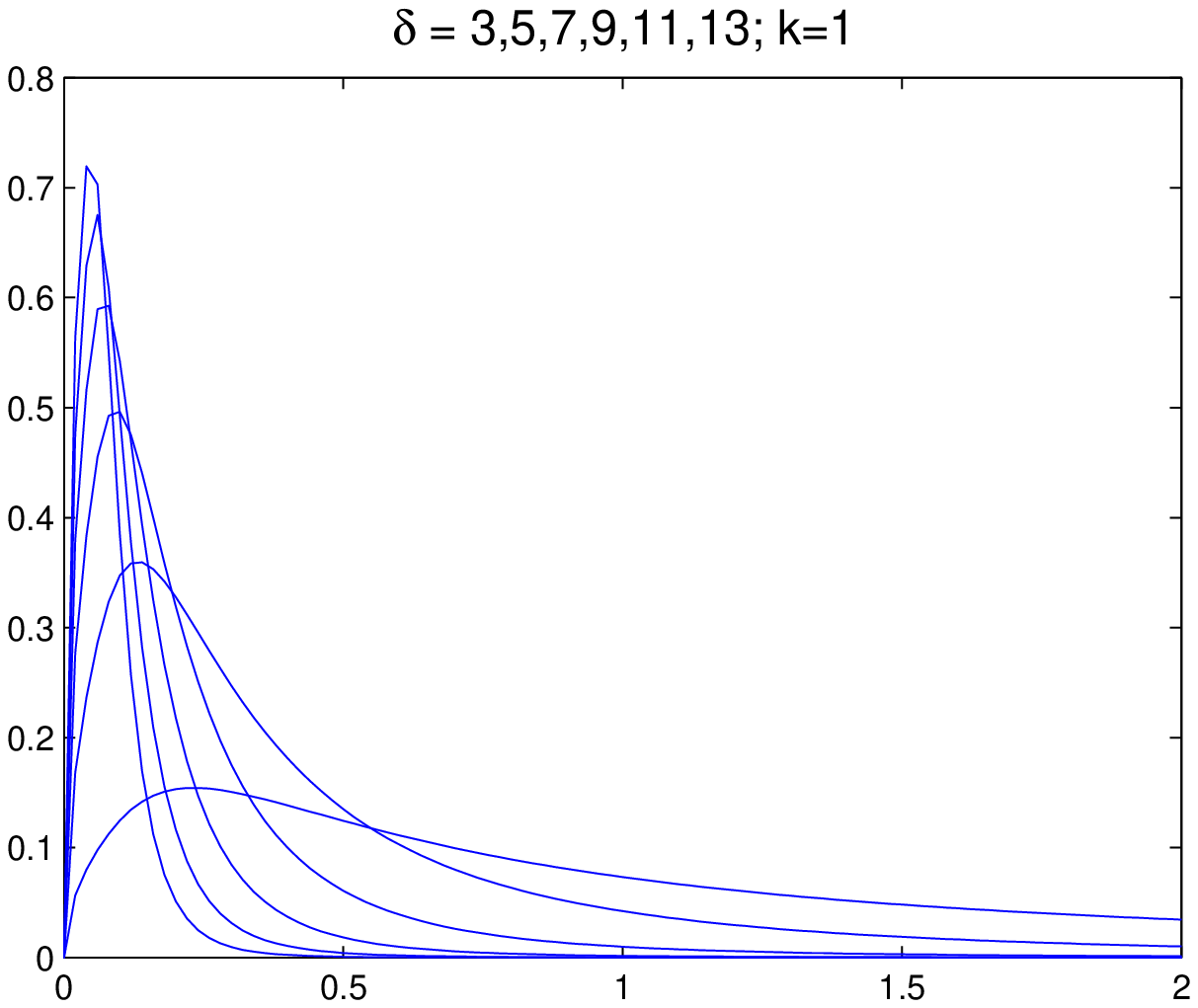,width=0.4\linewidth,clip=} &
\epsfig{file=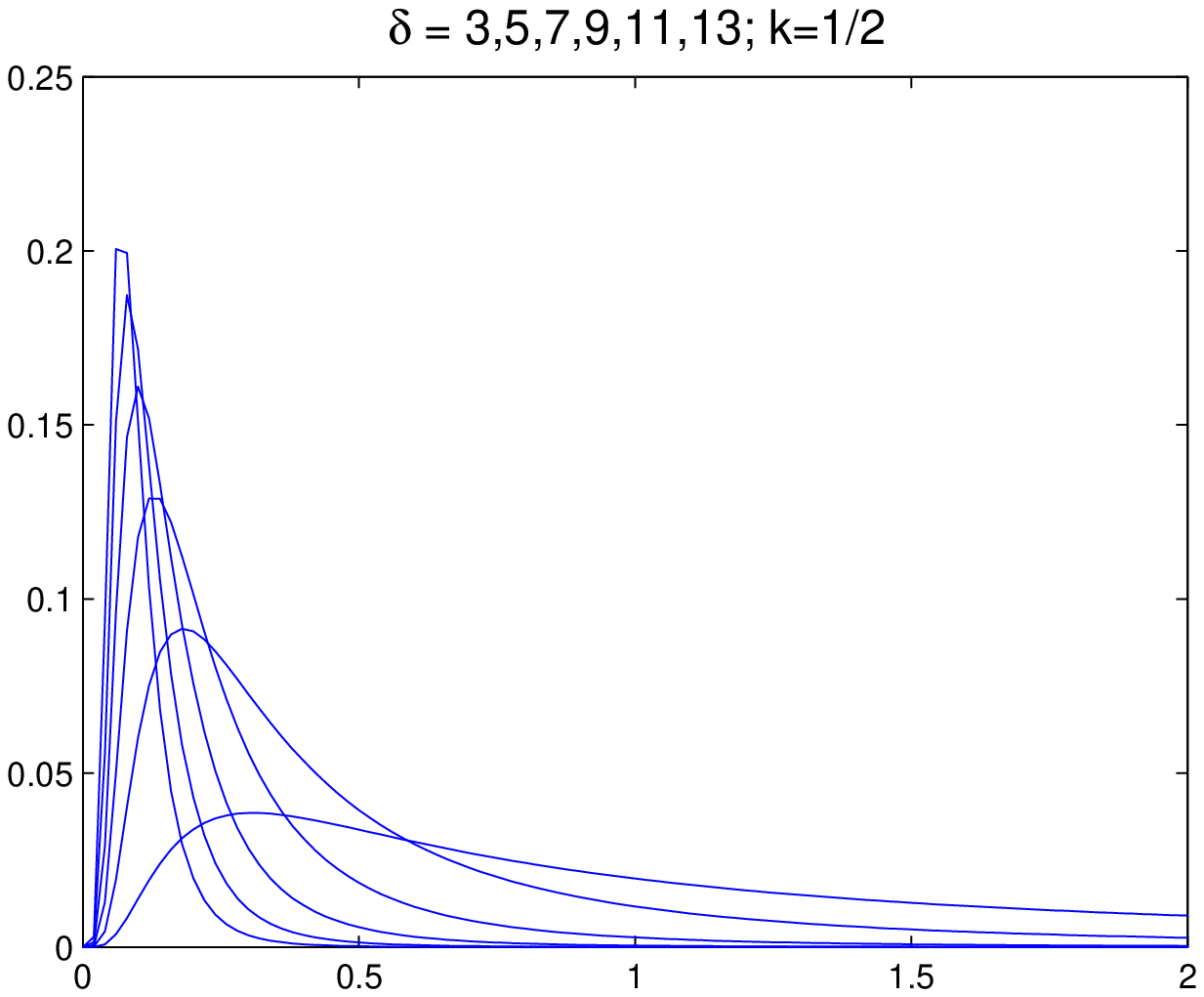,width=0.4\linewidth,clip=} \\
\epsfig{file=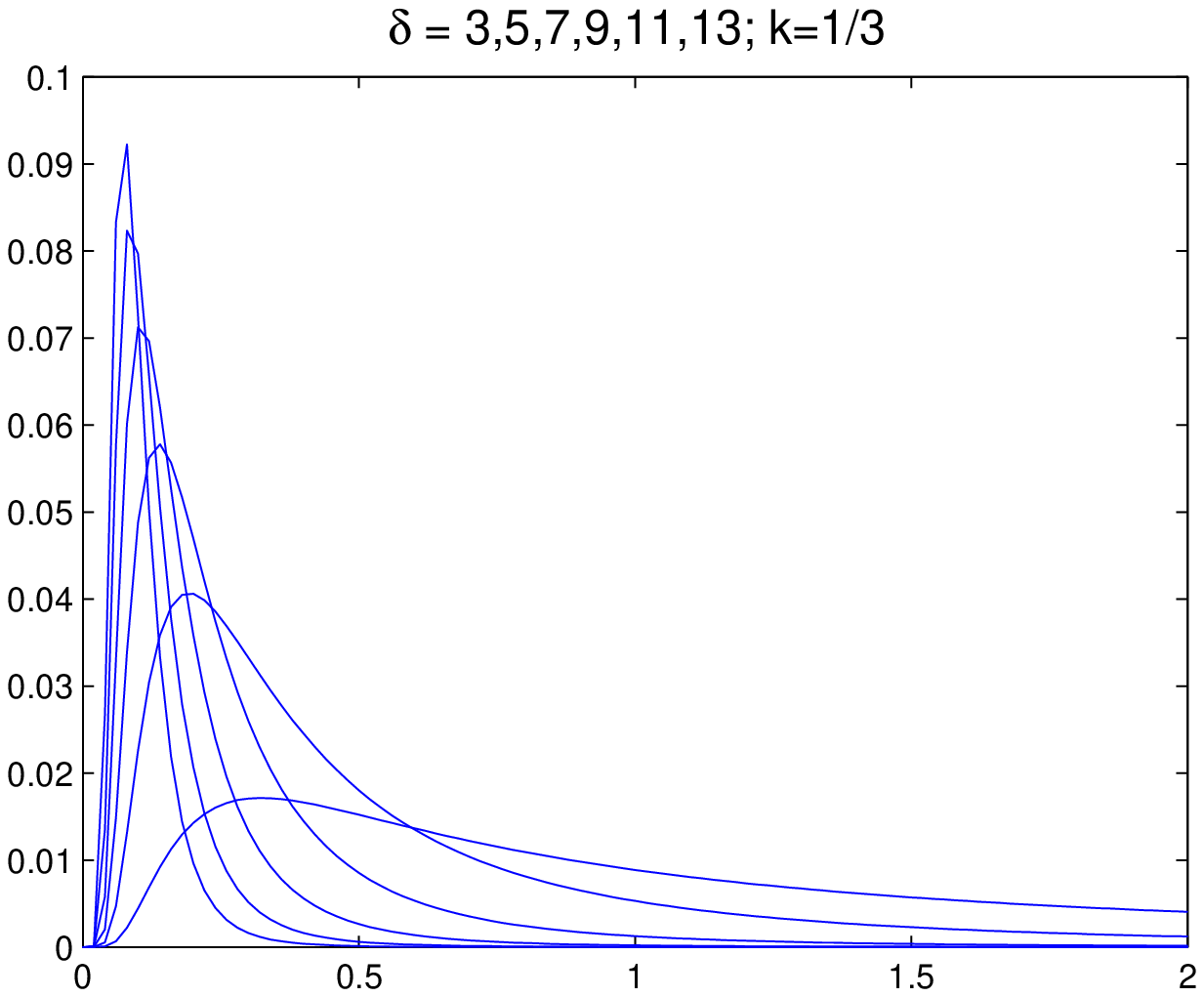,width=0.4\linewidth,clip=} &
\epsfig{file=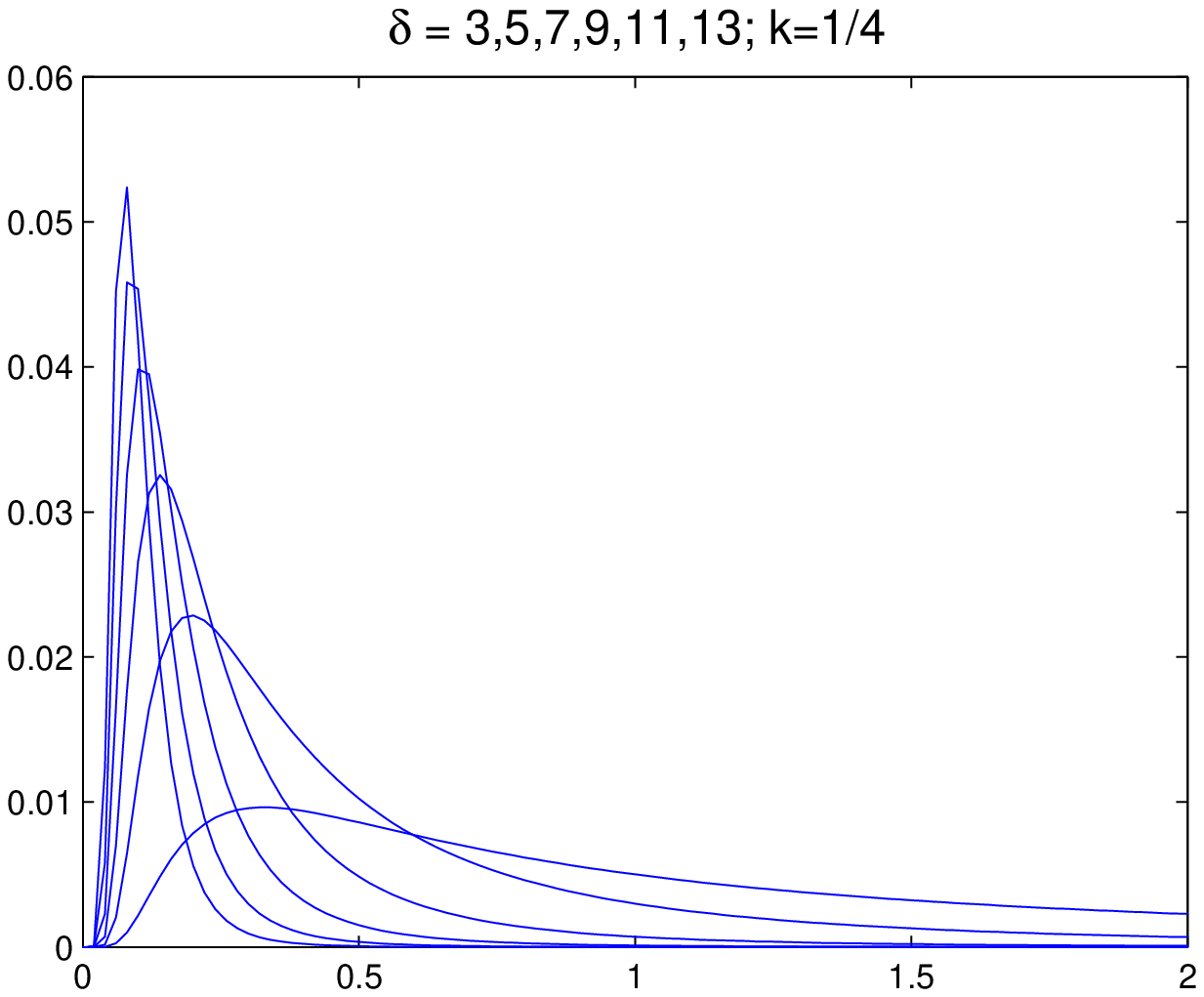,width=0.4\linewidth,clip=}\\
\epsfig{file=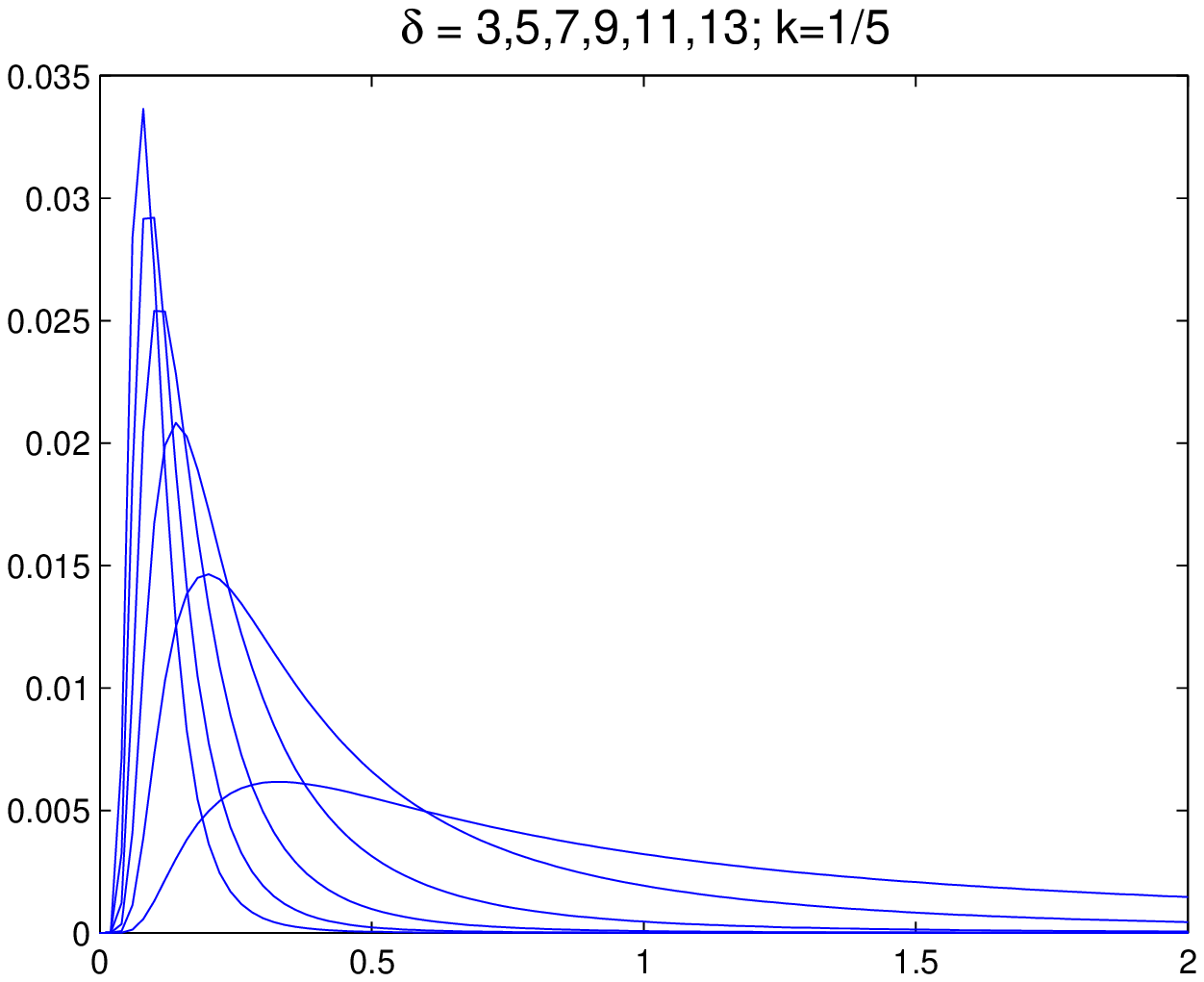,width=0.4\linewidth,clip=} &
\epsfig{file=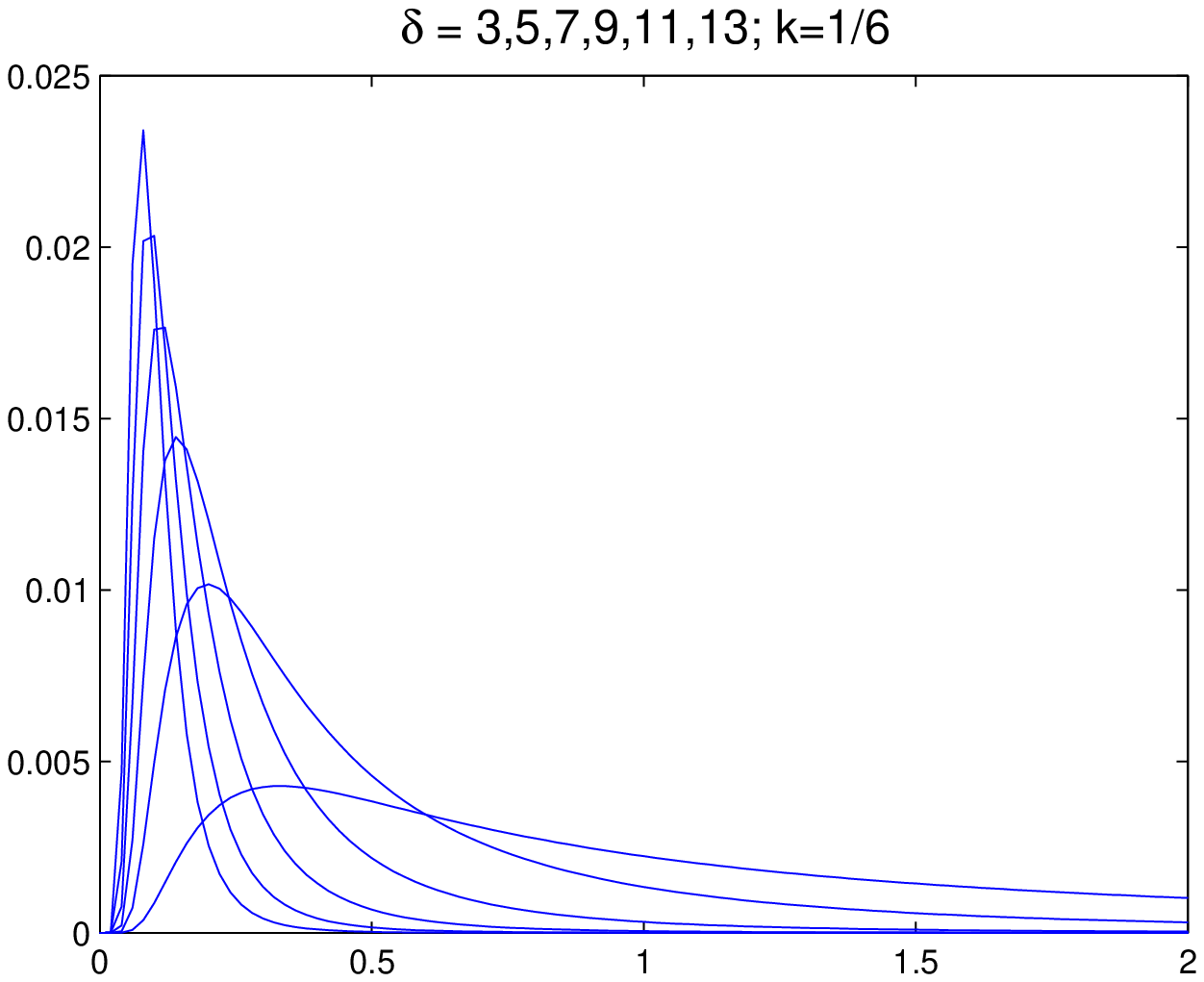,width=0.4\linewidth,clip=}\\
\epsfig{file=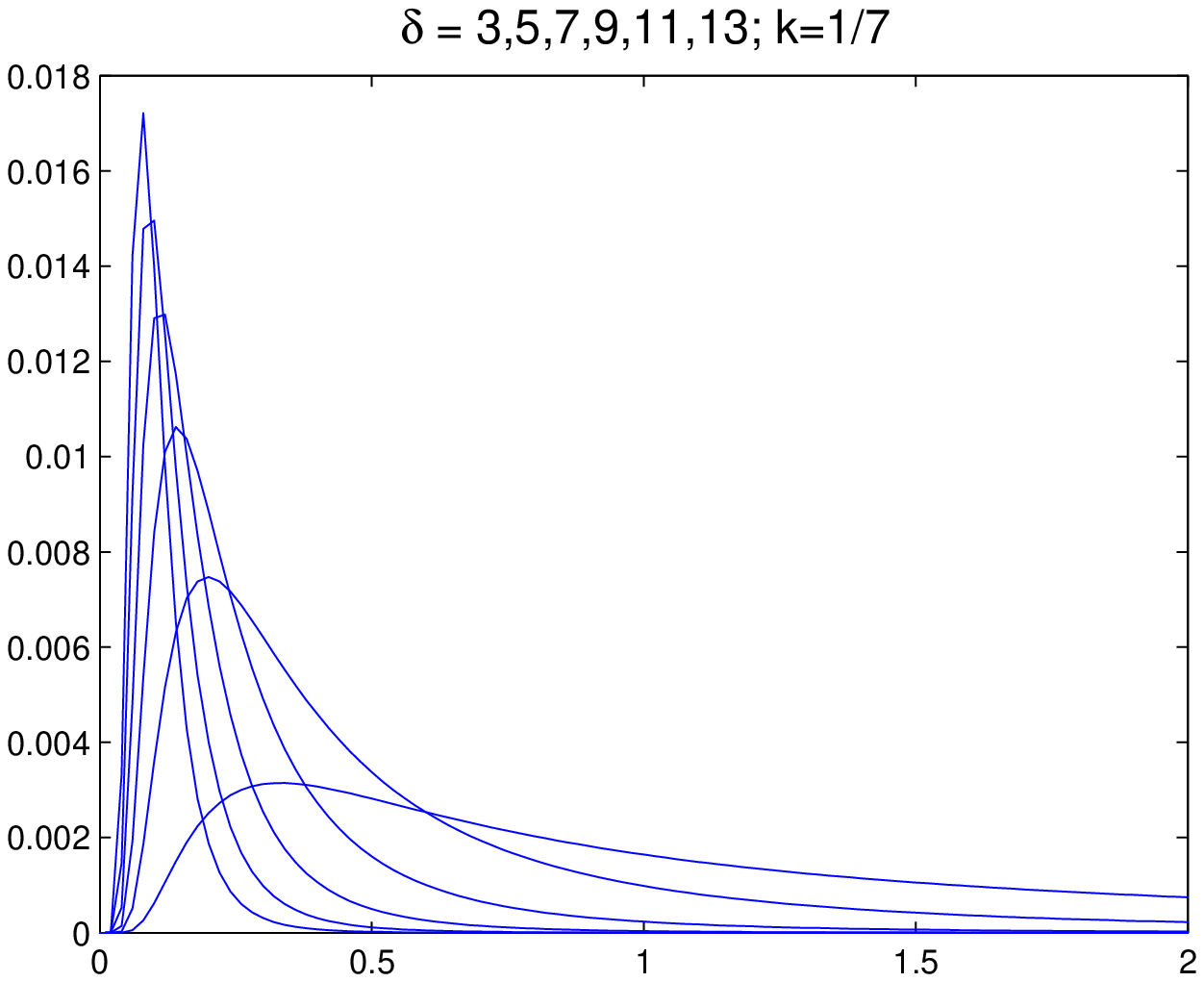,width=0.4\linewidth,clip=} &
\epsfig{file=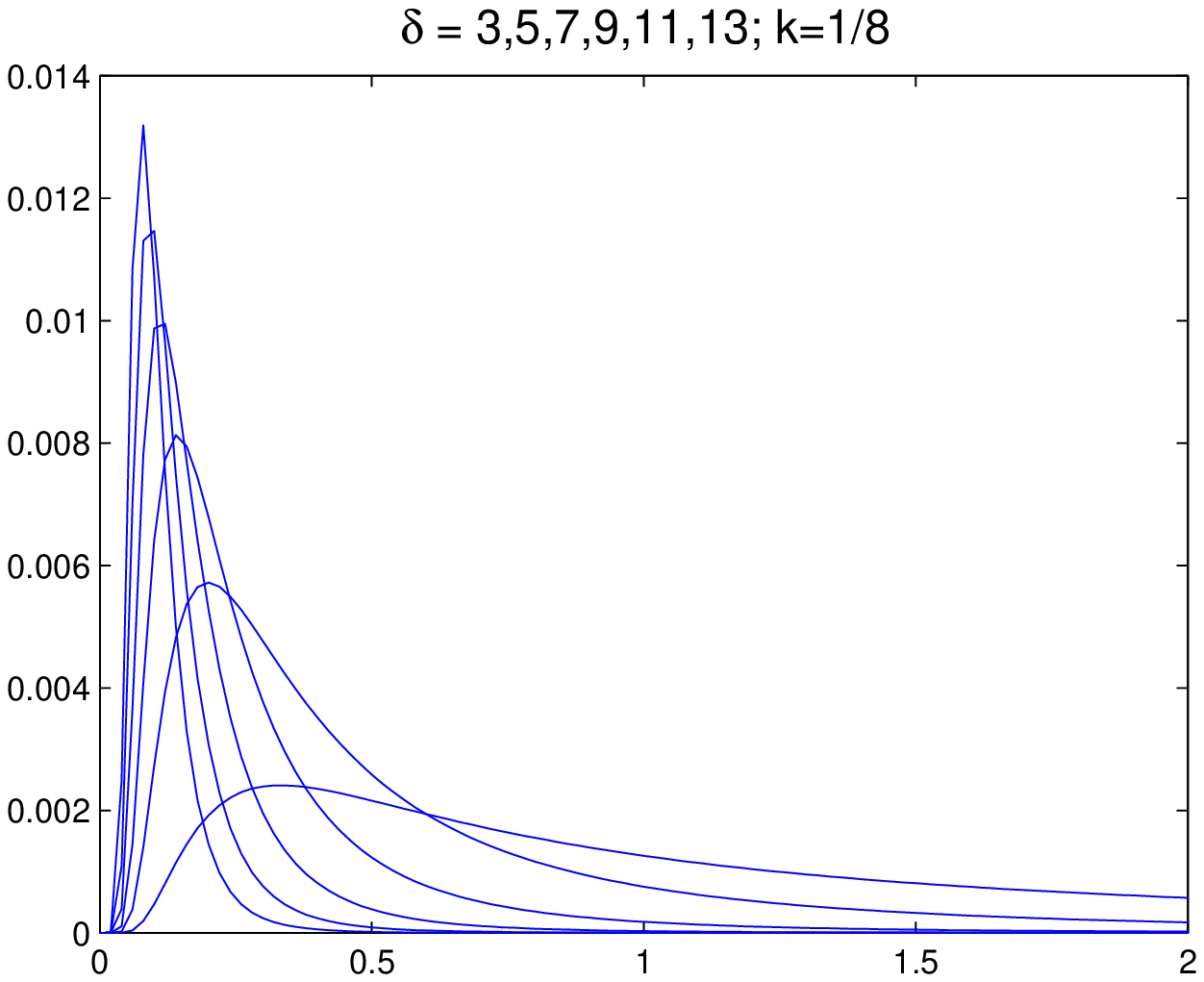,width=0.4\linewidth,clip=}\\
\epsfig{file=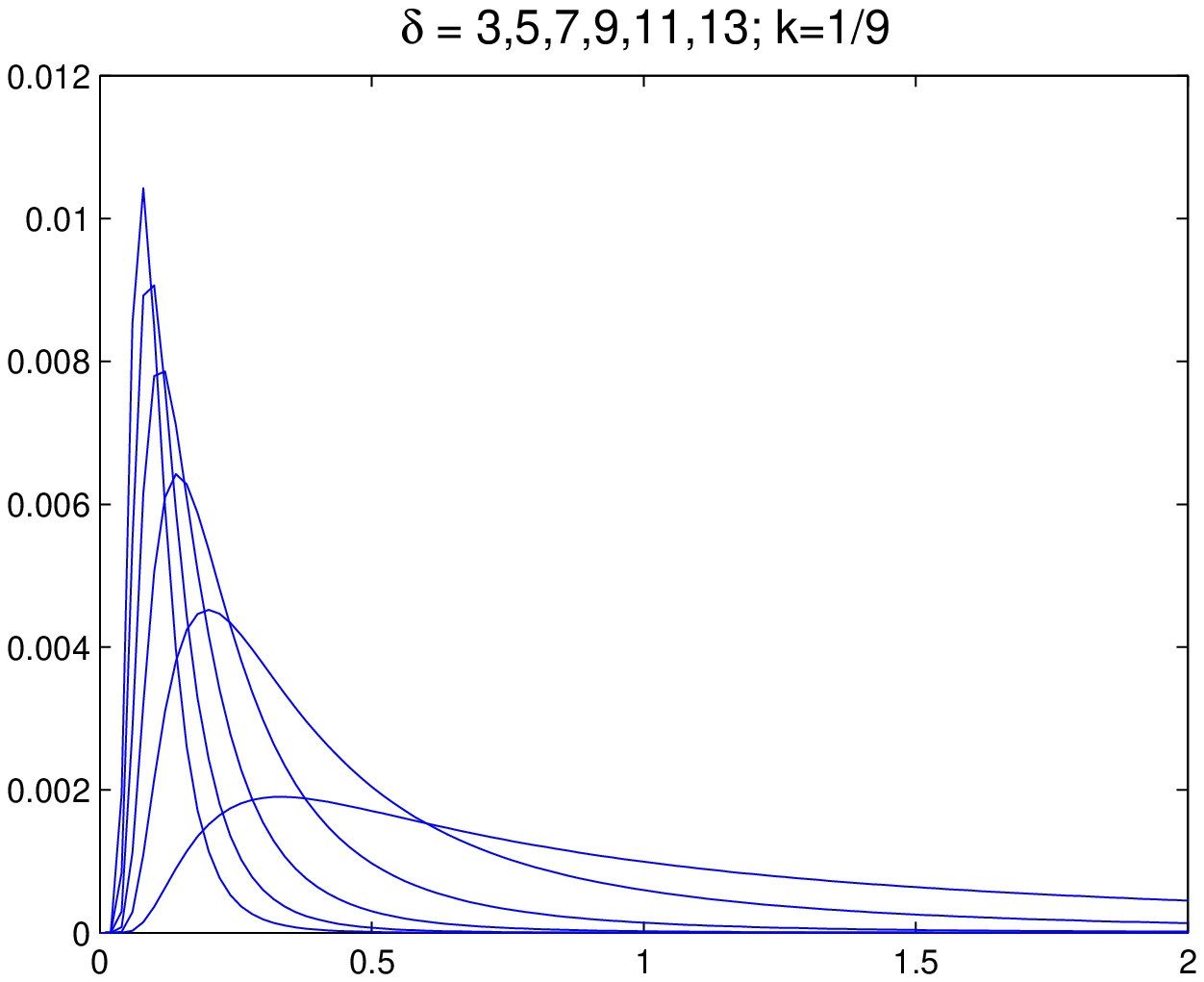,width=0.4\linewidth,clip=} &
\epsfig{file=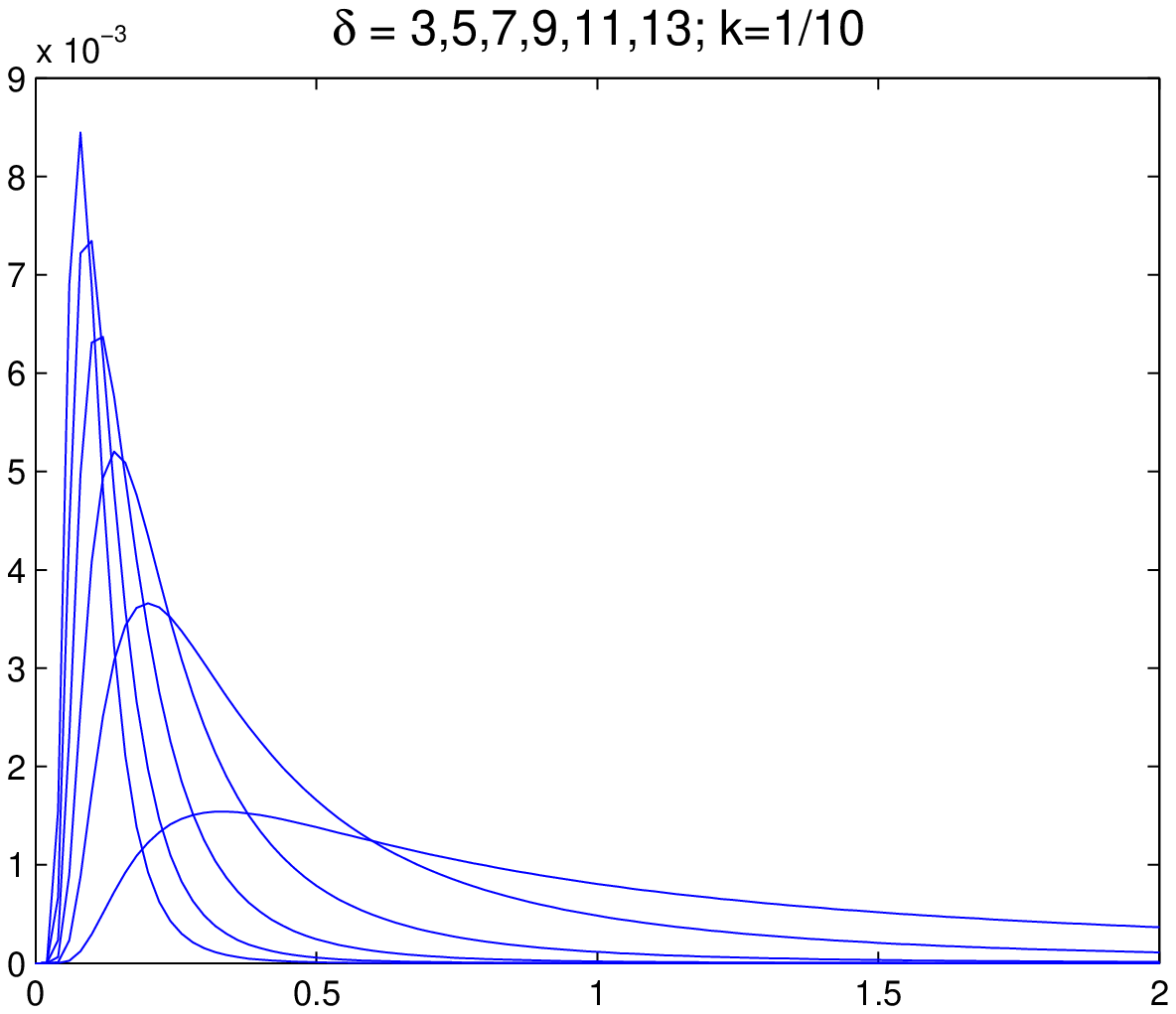,width=0.4\linewidth,clip=}
\end{tabular}
\caption{Varying $\delta$}
\end{figure}


\end{document}